\documentclass[a4paper,10pt]{article}
  \usepackage{amsmath,amsthm,amssymb, xcolor, amsfonts, epsfig, epstopdf, titling, url,array, enumerate, cite, mathrsfs, upgreek, authblk, scrextend}
  \usepackage[bindingoffset=0.2in,left=1in,right=1in,top=1in,bottom=1in,footskip=.25in]{geometry}
  \usepackage{graphicx}
  \usepackage{float, tikz}
  \theoremstyle{plain}

  \usepackage[bookmarksnumbered, colorlinks, plainpages]{hyperref}
  \hypersetup{colorlinks=true,linkcolor=violet, anchorcolor=green, citecolor=cyan, urlcolor=orange, filecolor=magenta, pdftoolbar=true}

  \newtheorem{thm}{Theorem}[section]
  \providecommand{\keywords}[1]{\begin{addmargin}[28pt]{28pt}\noindent\textbf{Keywords:} #1 \end{addmargin}}
  \newtheorem{lma}[thm]{Lemma}
  
  \newtheorem{ppn}[thm]{Proposition}
  \theoremstyle{definition}
  \newtheorem{dfn}[thm]{Definition}
  \newtheorem{eg}[thm]{Example}
  \newtheorem{rem}[thm]{\textit{Remark}}
  \providecommand{\ams}[1]{\begin{addmargin}[28pt]{28pt}\noindent\textbf{Mathematics Subject Classification:} #1\end{addmargin}}
  
  \title{Generalized $\theta$-Parametric Metric Spaces: Fixed Point Theorems and Applications to Fractional Economic Models}
  \author[1]{Abhishikta Das}
  \author[2]{Hemanta Kalita}
  \author[3]{Mohammad Sajid  \thanks{Corresponding author}}
  \author[4]{T. Bag}
  \affil[1,4]{ \small Department of Mathematics, Siksha-Bhavana, Visva-Bharati, Santiniketan-731235, West-Bengal, India }
  \affil[2]{ \small Mathematics Division, School of Advanced Sciences and Languages, VIT Bhopal University, Bhopal- Indore Highway, India }
  \affil[3]{ \small Department of Mechanical Engineering, College of Engineering, Qassim University, Saudi Arabia 
  	\authorcr 
  	Emails: abhishikta.math@gmail.com$^1$, hemanta30kalita@gmail.com$^2$, msajd@qu.edu.sa$^3$, tarapadavb@gmail.com$^4$
}
  \date{}
  
\begin{document}
  	\maketitle
  	
  	 \pagestyle{myheadings}
  	\markright{\footnotesize \it   Generalized $\theta$-Parametric Metric Spaces: Fixed Point Theorems and Applications to Fractional Economic Models}  
  	
  		\noindent\rule{\linewidth}{0.2pt}  
  	\begin{abstract}
  	The objective of this manuscript is to introduce and develop the concept of a generalized $\theta$-parametric metric space-a novel extension that enriches the   modern metric fixed point theory. We   study of its fundamental properties, including convergence and Cauchy sequences that establishes a solid theoretical foundation. 
  	A significant highlight of our work is the formulation   of   Suzuki-type fixed point theorem within this framework which extends   classical results in a meaningful way. 
  	To demonstrate the depth and applicability of our findings, we construct   non-trivial examples that illustrate the behavior of key concepts. Moreover, as a practical application, we apply our main theorem to analyze an economic growth model,    demonstrating its utility in solving fractional differential equations that arise in dynamic economic systems.\\
  	
  	\noindent
  	\ams{34A08; 34K37; 54E35; 47H10}
  \end{abstract}
  	\noindent\rule{\linewidth}{0.2pt}  

  	\section{Introduction}

  	The field of mathematics is developing at a very rapid pace these days. In particular, the study of derivatives and integrals, which were once limited to the natural number order, are now broadened to include the rational and real numbers in the fractional order. The field of fractional order integrals and derivatives is thought to have started in 1695, although major advances in the field have just been made in the early 21st century. Numerous scientific articles that study fractionally ordered derivatives and integrals and their applications in a variety of scientific domains, including   engineering and economics.   An FDE is often used to model growth incorporating memory effects. Since many economical processes have memory effects in their nature, the FDE is a suitable concept to model the growth of many economical processes.   Hence, calculus fractional serve as a backbone of the effect of memory on the economic growth model (please see the  references \cite{15, 16, 17, 18}).

  	Numerous areas of mathematical analysis, notably multidimensional calculus, real analysis, and complex analysis, are based on the general idea of metric spaces. The triangle inequality of the metric axiom is very fundamental and geometrically understandable and it makes us more comfortable with  metric spaces  than a general topological space.  That's why  the metric theory attracts  researchers.   As a result, by modifying or relaxing the  metric  axioms different kinds of  generalized metric spaces \cite{25}  have been developed. 
  	  To generalize the Banach contraction mapping theorem,   Czerwik \cite{5} provided an axiom that was less strong than the triangular inequality and officially defined a b-metric space.  As a generalization of metric spaces, Mustafa et al. \cite{3} created the structure of   G-metric spaces. 	A novel extension of the concept of metric space, which we refer to as an F-metric space, is presented by M. Jleli et al. in \cite{9}. They contrast their idea with preexisting literary generalizations. Additionally, in their scenario of $F$-metric spaces, a new form of the Banach contraction principle is constructed. In 2014, N. Hussain et al. \cite{11} defined parametric metric spaces where the distance function depends explicitly on a positive parameter. This concept provides additional flexibility in modeling various phenomena inspired by which A. Das et al. introduced   the idea of   generalized parametric metric in \cite{8}.

  	Another different kind of generalization of metric space   was done by F. Khojasteh et al. \cite{1} in 2013, where they developed the concept of $\theta$-metric space.   In this framework, they replaced the triangle inequality with a more generalized inequality using a family of binary operations on $ [ 0, \infty ) $, called  $ \mathcal{B} $-action.   They discussed  several important properties of the induced topology on such spaces, along with fixed point results for Banach-type and Caristi-type contraction mappings. The use of the  $\mathcal{B}$-action in the metric axiom has attracted considerable attention from researchers. Recently, Moussaouia, Radenovic  and Melliani \cite{12} introduced the concept of a $\theta$-parametric metric space. They elegantly applied the  	$\mathcal{B}$-action to the axiom of existing parametric metric spaces, thereby expanding the field of research on generalized metric spaces. 	With suitable examples, they demonstrated that the space they defined is fundamentally different from the parametric metric space. 
   This observation prompts us to explore the  appliance of the $\mathcal{B}$-action  to the axiom of generalized parametric metric spaces \cite{8}.  
   
    The outcomes derived from the uses of     $\mathcal{B}$-action in the definition of $\theta$-parametric metric spaces, when compared to the binary operation $`o$' discussed in \cite{8}, have inspired us to develop a more comprehensive framework, named as  `generalized $\theta$-parametric metric space', which may contribute on the expansion and  understanding the  metric space theories and their applications. 
   As demonstrated in the definition of a generalized parametric metric space, the modification of the role of the parameter significantly extends the axiom of the parametric metric space. Similarly, the generalized $\theta$-parametric metric space further extends the concept of the $\theta$-parametric metric space. This extension not only enriches the existing framework but also opens up new possibilities for exploring more complex structures within metric spaces.  
    As observed in the literature on generalized parametric metric spaces \cite{8,13}, the continuity of the binary operation $`o $' plays a crucial role in establishing fundamental results. In fact, every continuous binary operation qualifies as a $\mathcal{B}$-action. Henceforth, in the formulation of this new structure, we have opted to use a $\mathcal{B}$-action in place of a   binary operation $`o $', ensuring broader applicability and consistency with existing theoretical frameworks.
    In \cite{R1-9, R1-10, R1-11}, authors shown that the foundational study of this framework reveals an   fundamentally distinct approach compared to other generalized metric spaces. Inspired by its flexible structure, we adopt this concept in a novel manner within our present formulation.
    Our formulation aims to further explore the interaction between $\mathcal{B}$-action and parametric structures,   opening new avenues for future research in generalized metric spaces. Several recent contributions have significantly advanced fixed point theory in generalized metric frameworks, including  $\theta$-metric  and parametric metric spaces.  \cite{R1-1, R1-2, R1-3, R1-4,R1-5, R1-6, R1-7} are notable among them, which provide valuable insights and fixed point results in various extended metric contexts.

 	In the metric fixed point theory,  the Banach fixed-point principle has been widely studied.  In 2008, Suzuki \cite{10} presented the metric completeness of the underlying space using a fixed point theorem, which is a very basic generalization of the Banach contraction principle. Mathematicians called this fixed point theorem as the Suzuki fixed point theorem. The notion of an $F$-Suzuki contraction, or almost-$F$-contraction, was first presented by Wardowski \cite{24}, who also produced some intriguing fixed point results.  The broad applicability of $F$-Suzuki contractions motivates us to explore analogous Suzuki-type contraction principles within the framework of generalized $\theta$-parametric    metric spaces.   

  	This article is structured into four sections. Section 2 provides essential definitions and  results that will be helpful in the subsequent discussions. Section 3 introduces the notion of  a generalized metric, termed as generalized $ \theta $-metric space  which assigns a specialized class of functions to extend the triangle inequality. We compare this new framework with the generalized parametric metric  and $ \theta $-metric space, highlighting their distinctions and connections. Several illustrative examples are presented to support the proposed definition. 
  	Moreover, we establish conditions under which a generalized parametric metric can be induced from this new formulation.  	Additionally, we conduct a detailed examination of the  convergence of sequence and Cauchyness, continuity properties of the distance function, including topological structure of the space, such as openness and closedness of balls, compactness, sequential compactness.  We provide counterexamples to illustrate the limitations of the general setting and establish sufficient conditions under which classical topological properties are recovered. 
  	In Section 4, we develop a new Suzuki-type contraction principle within this framework and establish a corresponding fixed point theorem. As a significant application, we employ this theorem to investigate the existence of unique solutions to a Caputo fractional differential equation. 
  	The motivation for using Caputo fractional derivatives   stems from their ability to generalize classical differential equations while preserving essential properties such as initial conditions in the classical sense. In the framework of generalized $\theta$-parametric metric space, our aim to study over   fractional differential equations, particularly in the context of economic growth models.  This is not only highlights the theoretical significance of our results but also bridges the gap between abstract mathematical structures and real-world applications.

  	\section{Preliminaries} 
  	
  	In this Section, we recall several existing definitions and results. Throughout the article $ \mathcal{X}$ denote a non empty set.   We start the Section with the  definition of a b-metric space, which extends the triangle inequality through a constant multiplicative factor. 
  	
  	\begin{dfn} \cite{5}
  		A function $ b : \mathcal{X}^2 \to \mathbb{R}_{\geq 0 }$ is said to be an $b$-metric if and only if for all $x, \omega, \mu \in \mathcal{X}$ and all $r>0$ the following conditions are satisfies:
  		\begin{enumerate}[($b1$)]
  			\item $ b (x, \omega )=0$ if and only if $x= \omega;$
  			\item $ b (x, \omega )= b ( \omega ,x);$
  			\item $ b ( \omega, \mu ) \leq K[ b(x, \omega ) + b(x, \mu ) ] $ where $ K \geq 1 $ is a real constant.
  		\end{enumerate}
  		The pair $(\mathcal{X}, b )$ is called $b$-metric space.
  	\end{dfn}

   	Another  one  flexible structure  is $\theta$-metric space which based on the notion of a binary operation on non-negative real numbers, called $\mathcal{B}$-actions. The definition given below formalizes this concept.

  	\begin{dfn} \cite[Definition 4]{1} \label{def2}
  		Let $ \theta : [ 0, \infty )^2   \rightarrow [0, \infty ) $ be a continuous mapping with respect to each variables and $ Im ( \theta )  = \{ \theta ( a, \omega ) : a, \omega \in [ 0, \infty ) \} $. Then  $ \theta $ is called an   $ \mathcal{B}$-action  if $ \theta $  satisfies the followings:
  		\begin{enumerate}[($\mathcal{B}1$)]
  			\item $\theta (0, 0 ) = 0 $ and $\theta (x,  \omega ) = \theta (\omega, x ) $ for all $ x, \omega \geq 0 $;
  			\item $ \theta(x, \omega_1 ) <\theta(u, \omega_2 ) $ if either 		 $ x < \omega_2, ~   u \leq \omega_2   $ or $  x \leq \omega_1, ~   u < \omega_2    $; 
  			\item for each $ \eth \in  Im( \theta ) $ and for each $ x \in [0, \eth ] $, there exists $ \omega \in [0, \eth ] $ such that $ \theta (x, \omega ) = \eth $;
  			\item $ \theta ( \omega, 0 )  \leq  \omega $ for all $ \omega > 0 $.
  		\end{enumerate}
  	\end{dfn}

  	Several examples of  $ \mathcal{B}$-action   are listed below.

  	\begin{eg} \cite[Example 5]{1} \label{eg 2}
  		Let $ \omega, y \in [ 0, \infty ) $. Then following functions are some examples of  $\mathcal{B}$-action:  
  		\begin{enumerate}[(i)]
  			\item $\theta_1 (\omega, y ) = \omega + y + \omega y $,
  			\item   $ \theta_2 ( \omega, y ) = \dfrac{\omega y}{1 + \omega y } $,
  			\item  $\theta_3 ( \omega, y ) = k ( \omega + y + \omega y ) $, where  $ k \in ( 0, 1 ] $,
  			\item  $\theta_4 (\omega, y ) = ( \omega^2 + y^2 )^\frac{1}{2}$.
  		\end{enumerate}
  	\end{eg}

   Utilizing a  $\mathcal{B}$-action,  a new type of metric structure, known as a $\theta$-metric space \cite{1}, which replaces the classical triangle inequality with a $\mathcal{B}$-action, is  formally defined as follows:

  	\begin{dfn}\cite[Definition 11]{1} \label{dfn 3*}
  		Let $ \mathcal{X} $ be a nonempty set. A mapping $ d_\theta  : \mathcal{X}^2 \to  [0, \infty ) $	  is called a $ \theta $-metric  with respect to  $ \mathcal{B}$-action $ \theta  $ if   the followings hold:
  		\begin{enumerate}[($ d_\theta 1) $]
  			\item $ d_\theta (x, \omega ) = 0 $ if and only if $ x = \omega  $;
  			\item  $ d_\theta (x, \omega  ) =   d_\theta (\omega , x) $ for all $ x, \omega  \in \mathcal{X} $;
  			\item $ d_\theta (x, \omega  )  \leq  \theta( d_\theta (x, z ), d_\theta(z, \omega  ) ) $   for all $ x, \omega , z \in  \mathcal{X} $.
  		\end{enumerate}
  		A pair $ ( \mathcal{X}, d_\theta ) $ is called a $ \theta $-metric space
  	\end{dfn}

Next we recall the notion of a  parametric metric space, where the distance depends explicitly on a positive parameter.

\begin{dfn} \cite[Definition 5]{11}
  A  pair $( X, P)$ is called a parametric metric space   if  $ P : X \times X \times  (0,\infty) \rightarrow \mathbb{R}^+ $ be a function which 	 satisfies the following conditions:
	 \begin{enumerate}[(i)]
	 	\item $ P(a,b,t) = 0, ~ \forall~ t > 0 $ if and only if $ a = b $;
	 	\item $ P(a,b,t) = P(b,a,t), ~ \forall ~ a,b  \in  X $ and $ \forall~ t > 0 $;
	 	\item $ P(a,b,t) \leq P(a,x,t) + P(b,x,t), ~ \forall ~ a,b,x \in  X $ and $ \forall~ t > 0 $. 
	 \end{enumerate} 
\end{dfn}

  Recently, Moussaoui et al. \cite{12}  introduced the notion of    $\theta$-parametric metric by elegantly merging the axioms of the parametric metric \cite{11} with those of the $\theta$-metric \cite{1}. This notion unifies and extends these frameworks by incorporating a parameterized distance function governed by a $\mathcal{B}$-action, as described below.
  	
  	\begin{dfn} \cite[Definition 2.1]{12} \label{dfn 3}
  		A mapping $ d_\theta  : \mathcal{X}^2 \times (0, \infty ) \rightarrow  [0, \infty ) $ is called $ \theta$-parametric metric  on $  \mathcal{X} $ if
  		\begin{enumerate}[$(d_\theta 1)$]
  			\item $ d_\theta (\omega, y, p ) = 0  $ for all $ p >  0 $ if and only if $ \omega=y $;
  			\item $ d_\theta (\omega, y, p ) =  d_\theta ( y, \omega, p) $ for all $ p > 0 $ and  for all $ \omega, y \in \mathcal{X}  $;
  			\item $ d_\theta (\omega, y, p )  \leq  \theta \big( d_\theta (\omega,z,p), d_\theta(z,y,p) \big) $ for all $ p > 0 $ and for all $ \omega, y, z \in  \mathcal{X} $.
  		\end{enumerate}
  		Then the pair $ ( \mathcal{X}, d_\theta )$ is called a $\theta$-parametric metric space.
  	\end{dfn}

  \begin{rem}
  	Moussaoui et al. \cite{12} established that $\theta$-parametric metric in  general not a parametric metric but    every parametric metric  is effectively a $\theta$-parametric metric space with $\theta (u, v) =	u + v $.
  \end{rem}

\begin{eg} \cite[Example 2.3]{12}
  Let $ \mathcal{X} = \{x, y, z\} $ and  $\theta $ be a   $ \mathcal{B}$-action  defined by $  \theta (u, v) = u + v +  \sqrt{uv} $. A mapping $ \mathcal{D}_\theta  : \mathcal{X} \times \mathcal{X} \times (0, \infty ) \to [0, \infty ) $  defined by  
$$ \mathcal{D}_\theta (a, a, t ) =   0~ \forall a \in \mathcal{X}, ~ \mathcal{D}_\theta (x, y, t ) = \frac{1}{t}, \mathcal{D}_\theta (x, z, t) = \frac{3}{t}, ~ \mathcal{D}_\theta(z, y, t) = \dfrac{4 +  \sqrt{2}}{t}, $$ 
$$ \mathcal{D}_\theta(x, y, t) = \mathcal{D}_\theta(y, x, t), \mathcal{D}_\theta(x, z, t) = \mathcal{D}_\theta(z, x, t), \mathcal{D}_\theta(z, y, t) = \mathcal{D}_\theta (y, z, t) $$
for all $ t> 0 $. The pair $ ( \mathcal{X}, \mathcal{D}_\theta ) $ is $ \theta$-parametric metric space.	 
\end{eg}

 They also initiated the study of topological properties and the concept of convergence of  sequences by introducing the following definitions.

 \begin{dfn} \cite[Definition 2.5]{12}
   	Let $ ( \mathcal{X}, \mathcal{D}_\theta ) $ be a  $\theta $-parametric metric space. The open ball $ B_{\mathcal{D}_\theta} (x, r) $ with center $ x \in   \mathcal{X} $  and $ r \in Im(\theta ) $ 	is defined as:
  	$$	B_{\mathcal{D}_\theta}(x, r) = \{ y \in  \mathcal{X} :  \mathcal{D}_\theta (x, y, t ) < r \} ~ for ~ all ~ t > 0. $$ 
  \end{dfn}
  
  \begin{dfn} \cite[Definition 2.6]{12}
     Let $ ( \mathcal{X}, \mathcal{D}_\theta ) $   be a $\theta $-parametric metric space. A sequence $ \{ x_n \} $ in $ \mathcal{X} $ is said to be convergent to $ x \in  \mathcal{X} $, if for every $ \epsilon  > 0 $ there exists $ n_0 \in  \mathbb{N} $ such that $ \mathcal{D}_\theta (x_n, x, t ) < \epsilon $ for all $ n \geq n_0 $ and  $ t > 0 $.  
  \end{dfn}

\begin{dfn} \cite[Definition 2.7]{12}
 Let $ ( \mathcal{X}, \mathcal{D}_\theta ) $   be a $\theta $-parametric metric space.  
	\begin{enumerate}
		\item   A sequence $ \{ x_n \}$ in $ \mathcal{X} $ is called a Cauchy sequence if, $ \underset{ m, n \to \infty }{\lim} \mathcal{D}_\theta (x_n, x_m, t ) = 0 $ for all $ t > 0 $.
		\item  $ ( \mathcal{X}, \mathcal{D}_\theta ) $ is said to be complete if every Cauchy sequence $ \{ x_n \}$ in $ \mathcal{X} $ converges to some element in  $ \mathcal{X} $.
	\end{enumerate}
\end{dfn}

   The concept of generalized parametric metric spaces  rely on a special type of binary operations \cite{8} that generalize the additive structure  used in triangle inequalities. We begin by recalling the definition of such a binary operation.

  	\begin{dfn} \cite[Definition 3.1]{8} \label{dfn 4}
  		Let $ `o$' be a binary operation on $  [ 0,  \infty ) $  satisfying the followings:  
  		\begin{enumerate}[(a)]
  			\item $ z ~ o ~ 0 $ = $z$
  			\item $z \leq \beta \implies  z ~ o ~ a \leq \beta ~ o ~ a$ 
  			\item $z ~o ~ a  =  a ~ o ~ z$
  			\item $z ~ o ~( \beta ~ o ~ a ) = ( z ~ o ~ \beta ) ~o  ~a $
  		\end{enumerate}
  		for all $z, \beta, a  \in  [ 0,  \infty  ) $. 
  	\end{dfn}

 Based on such operations, Das et al.  \cite{8}  introduced the notion of a generalized parametric metric, as follows:

  	\begin{dfn} \cite[Definition 3.3]{8}\label{dfn 6}
  		Let $ \mathcal{X} $ be a non-empty set and $ \mathcal{P} : \mathcal{X}^2 \times  (0,  \infty )  \rightarrow   [0,  \infty ) $ be a function  satisfying
  		\begin{enumerate}[$ (\mathcal{P}{P}1) $]
  			\item $ ( \mathcal{P} (a, \nu, p) = 0, ~ \text{for all} ~ p > 0 )  $ if and only if $ a = \nu $;
  			\item $ \mathcal{P} (a,\nu,p) = \mathcal{P} (\nu,a,p)~ ~ \text{for all} ~  p > 0 $;
  			\item for $ p, t >  0 $,   $ \mathcal{P} (a,\nu,p+t) \leq \mathcal{P} (a, x, p) ~ o ~ \mathcal{P} (\nu,x,t) $;
  		\end{enumerate}
  		for all $ a,\nu, x \in  \mathcal{X} $.  
  		Then the function $ \mathcal{P} $ is said to be generalized parametric metric and $ ( \mathcal{X}, \mathcal{P}, o ) $ is called a generalized parametric metric space.
  	\end{dfn}

  	Recall several example of generalized parametric metric spaces below.

  	\begin{eg} \cite[Example 3.7]{8} \label{eg 3}
  		Consider two metrics $ \ss_{1}(a, \xi ) = \max \{ | a_{i} - \xi_{i} | : i= 1,2 \} $ and $ \ss_2(a, \xi ) = {\underset{i=1}{\overset{2} {\sum} } } |a_{i} - \xi_{i}|$ for all $ a = (a_{1}, a_{2} )$, $ \xi = ( \xi_{1}, \xi_{2} )  \in  \mathbb{R}^2$. 
  		Then the  function  $ \mathcal{P} $ defined by 
  		$$ \mathcal{P} (a, \xi, t) = \begin{cases}
  			100 ~~  \qquad & if \quad 0 < t \leq \ss_1(a, \xi ) \\
  			~	50 ~~  \qquad  & if \quad \ss_1(a, \xi ) < t \leq \ss_2 (a, \xi ) \\
  			~	25 ~~  \qquad &  if \quad   \ss_2(a, \xi ) < t \leq 2 \ss_2 (a, \xi ) \\
  			~~	0 ~~  \qquad  & if \quad 2 \ss_2 (a,\xi ) < t < \infty        
  		\end{cases} $$
  		for all $ a , \xi \in  \mathbb{R}^2 $ and for all $ t > 0 $, is a generalized parametric metric on $ \mathbb{R}^2 $ with respect to $`o$'= $`+$'. 
  	\end{eg}

  	\begin{eg} \cite[Example 3.7]{8}  \label{eg 4}
  		If $ ( \mathcal{X},d) $ be  a metric space then  the function $ \mathcal{P} :   \mathcal{X}^2 \times (0,\infty) \rightarrow [0,\infty)$  defined by  $ \mathcal{P} (a, \xi, p ) = e^p \ss ( a, \xi )$, for all $  a, \xi \in  \mathcal{X} $ and for all $ p > 0 $ is a parametric metric on $ \mathcal{X} $ but not a generalized parametric metric with respect to the binary operation $`+$'. 
  	\end{eg}

Some important properties of generalized parametric metric spaces are summarized below:  	
  	
  	\begin{ppn} \rm\cite[Proposition 3.9]{8}\rm \label{ppn 2}
  		If $ \mathcal{P} $ is a generalized parametric metric space on $ \mathcal{X} $. Then $ \mathcal{P} ( a, \xi, .) $ is non increasing function, for all $ a, \xi \in \mathcal{X} $. 
  	\end{ppn}

  	\begin{dfn} \cite[Definition 3.1]{13}  \label{dfn 8}
  		Let $ ( \mathcal{X}, \mathcal{P}, o ) $  be a generalized parametric metric space. For $ \iota > 0 $, the open ball $ B ( a, \alpha, \iota ) $ with center $ a \in \mathcal{X} $ and radius $ \alpha > 0 $ is defined by 
  		$	B(a, \alpha, \iota ) = \{\mu \in \mathcal{X} : \mathcal{P} (a, \mu, \iota )	< \alpha \}.$
  	\end{dfn}

  	\begin{thm} \rm\cite[Theorem 3.2]{13}\rm 
  		Let $ ( \mathcal{X},  \mathcal{P}, o ) $ be a generalized parametric metric space. Then   the collection of subsets of $\mathcal{X}$,  $ \tau_ {\mathcal{P}} = \{ \mathscr{O} \subseteq \mathcal{X} : ~ \text{for  any} ~ a \in \mathscr{O},  ~ \text{there exist} ~  \alpha, \iota \in ( 0, \infty)   ~ \text{such  that} ~ B(a, \alpha, \iota ) \subseteq \mathscr{O} \}  $   
  		is a  topology on  $ \mathcal{X} $.
  	\end{thm}

 There is another generalized metric space  called $\mathfrak{F}$-metric space which relies on a control function from a specific class of real-valued functions. We first define this class and then introduce the notion of an $\mathfrak{F}$-metric.

  	Let $\mathfrak{F}$ be the set of functions $ f : (0, \infty) \to \mathbb{R}$ satisfying the following conditions: 
  	\begin{enumerate}[($\mathfrak{F}_1$)]
  		\item $ f  $ is non-decreasing;
  		\item For every sequence $\{ \varepsilon_n\} \subset (0, \infty ) $, we have 
  		$$ \lim_{n \to +\infty} \varepsilon_n = 0 \iff \lim_{n \to +\infty} f( \varepsilon_n ) = -\infty. $$
  	\end{enumerate}
  	
  Throughout the article, we denote $\mathfrak{F}$ for mentioning the family of such functions.	
  	
  	\begin{dfn} \cite[Definition 2.1]{9} \label{dfn 9}
  		Let $ \mathcal{X} $ be a non-empty set and $ \ss : \mathcal{X}^2 \to [0, \infty ) $ be a given mapping. Suppose that there exists $ ( f, \alpha ) \in  \mathfrak{F} \times  [ 0, \infty ) $ such that 
  		\begin{enumerate}[$ (\ss1) $]
  			\item $ ( x, \eta )  \in  \mathcal{X} \times \mathcal{X}, ~ \ss ( x, \eta ) = 0 \iff x = \eta $;
  			\item $ \ss ( x, \eta ) = \ss ( \eta, x ) $ for all $ ( x, \eta ) \in  \mathcal{X} \times  \mathcal{X} $;
  			\item For every $ ( x, \eta ) \in  \mathcal{X} \times \mathcal{X} $, for every  $ \mathfrak{N} \in \mathbb{N} $, $ \mathfrak{N} \geq 2 $ and for every $\{ u_\iota \}_{ \iota = 1 }^\mathfrak{N} \subset \mathcal{X} $ with $(u_1, u_\mathfrak{N} ) = (x, \eta )$, we have 
  			$$ \ss ( x, \eta ) > 0 \implies f( \ss (x, \eta ) ) \leq f \left( {\underset{i=1}{\overset{ \mathfrak{N} - 1 }{ \sum } }} \ss (u_\iota, u_{ \iota + 1 }) \right) +\alpha. $$
  		\end{enumerate}
  		Then $ \ss $ is said to be an $\mathfrak{F}$-metric on $ \mathcal{X} $ and the pair $ (\mathcal{X}, \ss ) $ is said to be an $\mathfrak{F}$-metric space.
  	\end{dfn}

  As we make use of a Suzuki-type contraction to establish a  fixed point result in the generalized $\theta$-parametric metric space, before extending it to our setting, we recall its classical version in the framework of metric spaces, as introduced by Suzuki \cite{10}.

  	\begin{thm} \rm\cite[Theorem 2]{10}\rm
  		\label{suzuki in m s}
  		Let $ \mathfrak{T} $ be a self mapping on a complete metric space $ (\mathcal{X}, \ss ) $  and $ \psi : [ 0, 1) \to ( \frac{1}{2}, 1 ] $ be defined by 
  		$$ \psi ( \rho ) = \begin{cases}
  			~~ 1 \qquad ~~~~ if ~~~ 0 \leq \rho \leq \dfrac{\sqrt{5} - 1}{2} \\
  			\dfrac{1 - \rho}{\rho^2} \qquad if ~~~ \dfrac{\sqrt{5} - 1}{2} \leq \rho \leq \dfrac{1}{\sqrt{2}} \\
  			\dfrac{1}{\rho + 1} \qquad if ~~~ \dfrac{1}{\sqrt{2}} \leq \rho < 1.
  		\end{cases}
  		$$
  		If there exists $ \rho \in [ 0, 1) $ such that for all $ \mu, \eta \in \mathcal{X} $,
  		$$ \psi ( \rho ) \ss ( \mu, \mathfrak{T} \mu ) \leq \ss (\mu, \eta ) \implies \ss ( \mathfrak{T} \mu, \mathfrak{T} \eta ) \leq \rho \ss (\mu, \eta ) $$ 
  		then $ \mathfrak{T} $ has an unique fixed point in $\mathcal{X}$.
  	\end{thm}

  	\section{Introduction to generalized $ \theta $-parametric metric space} 
  	
  	In this Section of the article, we develop the concept of generalized $ \theta $-parametric metric space and study several characterization, beginning with its formal definition.

  	\begin{dfn}\label{dfn 1}
  		Let $ \mathcal{X} $ be a non empty set and $ \theta $ be a   $ \mathcal{B} $-action  on $ [ 0, \infty ) $ as defined in the  Definition \ref{def2}. Suppose $ \mathfrak{P}_{\theta}  : \mathcal{X} \times \mathcal{X} \times ( 0, \infty )  \rightarrow [ 0, \infty )  $ be a function which satisfies the  following conditions:
  		\begin{enumerate}[$(\mathfrak{P}_{\theta}1)$]
  			\item $ ( \mathfrak{P}_{\theta} (x, \omega, p ) = 0 ~ \forall p > 0 ) $ iff $ x = \omega $;
  			\item there exist $ ( f, \alpha )  \in \mathfrak{F} \times [0,\infty) $ such that for all $ s, p > 0 $ and $ x, \omega, \mu \in \mathcal{X} $, 
  			$$ ( \mathfrak{P}_{\theta} ( x , \omega, r ) > 0 ~ \forall r > 0 ) \implies  f ( \mathfrak{P}_{\theta} ( x, \omega, s + p ) ) \leq  f \left(    \theta   ( \mathfrak{P}_{\theta } ( x, \mu, s ), \mathfrak{P}_{\theta } ( \mu, \omega, p ) )    \right)  + \alpha. 
  			$$ 
  		\end{enumerate}	
  		We called $ \mathfrak{P}_{\theta} $ is a generalized $ \theta $-parametric metric   and  the  tuple $ (\mathcal{X}, \mathfrak{P}_{\theta}, \theta ) $ is called a generalized $ \theta $-parametric metric space with respect to the pair $ ( f, \alpha )  \in \mathfrak{F} \times [0,\infty) $.
  	\end{dfn}
  	Next we  present some examples  of generalized $ \theta $-parametric metric space as follows:

  	\begin{eg}
  		Let us define a function $ \ss : \mathbb{Z} \times \mathbb{Z} \to [0, \infty )  $ by 
  		$$ \ss ( x, \xi ) = \begin{cases}
  			~ | x - \xi |  \quad ~~~~~~~ \text{if} ~ x, \xi ~ \text{both are even or odd} \\
  			5 | x - \xi | + 9  \quad ~~ \text{if  one of}~ x, \xi ~ \text{is even or odd}.
  		\end{cases} $$
  		It is not hard to see $ \ss ( x, \xi ) \leq 5 ( \ss ( x, z ) + \ss ( z, \xi ) ) $ holds for all $ x, \xi, z \in \mathbb{Z} $. 
  		
  		Next we define a function $ \mathfrak{P}_{\theta} : \mathbb{Z} \times \mathbb{Z} \times ( 0, \infty ) \to [0, \infty )  $  by $ \mathfrak{P}_{\theta} ( x, \xi, q ) = \frac{\ss ( x, \xi ) }{q} $ for all $ x, \xi \in \mathbb{Z} $ and for all  $  q > 0 $.
  		Let us consider $ \mathcal{B}$-action as $ \theta ( \eth_1, \eth_2 ) = \eth_1 + \eth_2 $ whenever for all $ \eth_1, \eth_2 \in [0, \infty ) .$ 
  		
  		Then clearly $ \mathfrak{P}_{\theta} $ satisfies $ (\mathfrak{P}_{\theta} 1) $. 
  		For $ (\mathfrak{P}_{\theta} 2) $, let $ x, \xi, z \in \mathbb{Z} $ such that $ x \neq \xi $. By Definition,   for all  $ v > 0, ~ \mathfrak{P}_{\theta} ( x, \xi, v ) \neq 0 .$ Thus, for any   $ s, t > 0 $ we have 
  		\begin{align*}
  			&	\mathfrak{P}_{\theta} ( x,\xi, s+ t ) =   \frac{ \ss ( x, \xi ) }{ s + t } 
  			\leq    ~ \frac{ 5 ( \ss ( x, z ) + \ss ( z, \xi ) ) }{ s + t }   
  			<    ~ 5 \left( \frac{   \ss ( x, z )  }{ s  } +  \frac{    \ss ( z, \xi )  }{  t }  \right) \\
  			\implies & \ln ( \mathfrak{P}_{\theta} ( x,\xi, s+ t ) ) < \ln 5 + \ln \bigg(     \theta (	\mathfrak{P}_{\theta} ( x, z, s ) ,  	\mathfrak{P}_{\theta} ( z, \xi,  t ) )    \bigg).
  		\end{align*} 
  		
  		Therefore $ ( \mathbb{Z}, \mathfrak{P}_{\theta}, + ) $ is a generalized $ \theta $-parametric metric space with respect to $  f ( \mu ) = \ln \mu, ~ \mu > 0 $ and  $ \alpha = \ln 5 $.
  	\end{eg}

  	\begin{eg}
  		Let $ k  \in [1, \infty ) $  be a fixed positive real  number. We define a function  
  		$ \mathfrak{P}_{\theta} : \mathbb{R}^2 \times ( 0, \infty ) \to [ 0, \infty ) $ by 
  		$$ \mathfrak{P}_{\theta} ( x, \beta, t ) = \begin{cases}
  			k e^{\frac{| x- \beta |}{t}}  \quad if ~~ x \neq \beta \\
  			~~ 0 \quad ~~~~~~ if ~~  x = \beta 
  		\end{cases}
  		$$
  		for all $ x, \beta \in \mathbb{R} $ and $ t > 0 $. In this occasion, we consider $ \beta$-action by $ \theta ( \eth_1, \eth_2 ) = \max \{ \eth_1, \eth_2 \} $ for all $ \eth_1, \eth_2 \in [ 0, \infty ) $. Then for all $ x, \beta, z \in \mathbb{R} $ with $ x \neq \beta $ and $ s, t > 0 $, we have
     	\begin{align*}
  			\frac{1}{k} - \frac{1}{ \theta ( \mathfrak{P}_{\theta} ( x, z, s ) ,  	\mathfrak{P}_{\theta} ( z, \beta,  t )  )} + \frac{1}{ \mathfrak{P}_{\theta} ( x, \beta, s+ t )}
  			= & \frac{1}{k} - \frac{1}{ \max \bigg\{ k e^{\frac{| x- z |}{s}}, k e^{\frac{| z - \beta |}{t}}\bigg\} } +  \frac{1}{ k e^{\frac{| x- \beta |}{ s+ t}} }   \\
  			> & \frac{1}{k} - \frac{1}{ \max \bigg\{ k e^{\frac{| x- z |}{s}}, k e^{\frac{| z - \beta |}{t}}\bigg\} }. 
  		\end{align*}
  		Since $  \max \bigg\{ k e^{\frac{| x- z |}{s}}, k e^{\frac{| z - \beta |}{t}} \bigg\} \geq  k e^{\frac{| x- z |}{s}} \geq k  $, above relation gives 
  		\begin{equation}\label{new1}
  			\frac{1}{k} - \frac{1}{ \theta ( \mathfrak{P}_{\theta} ( x, z, s ) ,  	\mathfrak{P}_{\theta} ( z, \beta,  t )  )} + \frac{1}{ \mathfrak{P}_{\theta} ( x, \beta, s+ t )}  \geq   0.
  		\end{equation} 
  		Again by definition of $ \mathfrak{P}_{\theta} $, $ \mathfrak{P}_{\theta} ( x, \beta, s + t  ) > 0 $. Next using the function  $ f (\eta ) = - \frac{1}{\eta }, ~ \eta > 0 $ in (\ref{new1}), we get	
  		$$ f ( \mathfrak{P}_{\theta} ( x, \beta, s+ t ) ) \leq f ( \theta ( \mathfrak{P}_{\theta} ( x, z, s ),  	\mathfrak{P}_{\theta} ( z, \beta,  t )  ) ) + \frac{1}{k}.
  		$$
  		Hence $ ( \mathbb{R}, \mathfrak{P}_{\theta}, \max )     $ is a generalized $ \theta $-parametric metric space with respect to $ f ( \eta ) = - \frac{1}{\eta }, ~ \eta > 0 $ and $ \alpha = \frac{1}{k} .$ 
  	\end{eg}

  	\begin{rem}  
  	   A generalized $ \theta $-parametric metric space  reminisces  applience of parameter  as in fuzzy metric space. Although both fuzzy and parametric metric spaces aim to generalize the classical notion of a metric by introducing a parameter $t \in (0, \infty)$, their underlying axiomatic structures are fundamentally distinct. 
  		A fuzzy metric space, in the sense of George–Veeramani \cite{R1-8}, is  a triplet $(X, M, \star)$, where $M(x, y, t) \in [0,1]$ for all $x, y \in X$ and $t > 0$,  represents the grade of the  metric function $ d $ of $ x $ and $ y $ w.r.t. $ t $ so that $ d(x,y) $ less than $ t $ holds,  governed by a continuous t-norm $`\star$'. This framework is designed to model uncertainty and imprecision in measurements.
  		 In contrast, a parametric metric space $(X, P)$ \cite{11} also involves a parameter $t > 0$, but the function $P(x, y, t)$ is real-valued and deterministic. It generalizes the classical metric by embedding the metric behavior into a parametric structure, without incorporating fuzziness or probabilistic interpretation.
  		
    Therefore, despite their common dependability on a parametric approach, there is no direct internal dependency between fuzzy metric and parametric metric. This distinction highlights that the generalized $\theta$-parametric metric framework offers methodologies and applications that are separate and independent from those found in fuzzy metric space theory.
  \end{rem}

  	\begin{rem} \label{rem 1}
  		If the given binary operation $`o$' is  continuous and  the  $ \mathcal{B} $-action  is being $ \theta  $, then  every generalized parametric metric space is a generalized $ \theta $-parametric metric space (see \cite{8}).  For justification, let $ ( \mathcal{X}, \mathcal{P}, o )$ be a generalized parametric metric space. Clearly $\mathcal{P}$ is a non-negative real valued function satisfying the axiom $ ( \mathfrak{P}_{\theta}1)$. 
  		Again for all $ \xi, y, \mu \in \mathcal{X} $ with $  \xi \neq y $ and $ s,t > 0 $, the  parametric  metric   $ \mathcal{P} $ satisfies
  		$$ \mathcal{P} ( \xi, y, s+t) \leq \mathcal{P} ( \xi, z, s) ~ o ~  \mathcal{P} ( \mu, y, t) 
  		\implies \ln \mathcal{P} (\xi, y, s+t) \leq \ln (\theta ( \mathcal{P} (\xi, \mu, s), \mathcal{P} ( \mu, y, t ) ). $$
  		Therefore $ \mathcal{P} $ satisfies $ (\mathfrak{P}_{\theta}2)$  for $ f( \mu ) =\ln \mu, ~  \mu > 0 $ and $ \alpha = 0 $.
  	\end{rem}

  	It is worthy to mention that $ \theta $-parametric metric space and generalized $ \theta $-parametric metric space are different concept of distance functions.  We justify that there is a generalized $ \theta $-parametric metric space which is not a  $ \theta $-parametric metric space in the following example.

  	\begin{eg}\label{eg 1} 
  		Consider the generalized parametric space $ ( \mathcal{X}, \mathcal{P}, o ) $ of the Example \ref{eg 3}. Then $ \mathcal{P} $ is   a generalized $ \theta $-parametric metric space with respect to $ f(t) =\ln t $  and  $ \alpha =0 $ under the $ \mathcal{B}$-action $ \theta ( \eth_1, \eth_2 ) = \eth_1 + \eth_2 $ for all $ \eth_1, \eth_2 \in [0, \infty ) $. On the other hand, for $  {a}=(1,0), ~ \sigma = ( 0, \dfrac{1}{2}) $ and $ \tau = ( \dfrac{1}{4} ,\dfrac{1}{8} ) $ and $ t=1 $, one can see $	\theta ( \mathcal{P} ( a, \tau, t ), \mathcal{P}( \sigma, \tau, t ) ) 	=    \mathcal{P} ( a, \tau, t ) + \mathcal{P} ( \sigma, \tau, t )   =  50 < 100  =   \mathcal{P} ( a, \sigma, t ). $ Therefore $  \mathcal{P} $ is not a  $\theta $-parametric metric on $ \mathcal{X} $.
  	\end{eg}
  	
  	Further, a generalized  parametric metric $ \mathcal{P} ( x, \mu, t) $ on $ \mathcal{X} $   is, in a general sense, an  increasing function of $ t > 0 $, for all $ x, \mu \in \mathcal{X} $. However, a generalized $\theta $-parametric metric does not necessarily retain this property. To illustrate this, we present the following example.

  	\begin{eg}\label{inc-dec p}
  		Suppose $ ( \mathcal{X}, \ss )$ be a metric space, and  $ \mathfrak{P}_{\theta} : \mathcal{X}^2 \times ( 0, \infty ) \to [0, \infty )  $  is defined  by
  		$$ \mathfrak{P}_{\theta} ( a, \sigma, t ) = \begin{cases}
  			~~~ 25 ~~ \hskip 40 pt \mbox{if} ~~~  0 < t \leq \ss ( a, \sigma ) \\
  			~~~ 50  ~~ \hskip 40 pt \mbox{if} ~~ ~ \ss ( a, \sigma ) < t \leq 2 \ss ( a, \sigma )   \\
  			\dfrac{100 d ( a, \sigma ) }{t} ~ ~~~~ ~~ \mbox{if} ~~~ 2 \ss ( a, \sigma ) \leq t < \infty
  		\end{cases} $$
  		for all $ a, \sigma \in \mathcal{X} $ and $ t > 0 $.  We shall establish  $ \mathfrak{P}_{\theta} $ is a generalized  $\theta $-parametric metric on $\mathcal{X}$ under $ \mathcal{B} $-action  $ \theta ( \eth_1, \eth_2 ) = \dfrac{ \eth_1 + \eth_2 }{ 2 }  $ for all $ \eth_1, \eth_2 \in [0, \infty )  $. 
  		It is easy to see $  \mathfrak{P}_{\theta} $ satisfies the axiom $ ( \mathfrak{P}_{\theta} 1)	.$ We shall only verify the axiom $ ( \mathfrak{P}_{\theta} 2)	$. Let $ a, b, c \in \mathcal{X} $ so that $ a \neq b $. Then for all $ v > 0 $, $ \mathfrak{P}_{\theta} ( a, b, v ) > 0 $. Now choose   $ s, t > 0 $. Then we have the following cases:
  		\begin{enumerate}[(i)]
  			\item Let $  2 \ss ( a, c ) \leq s < \infty, ~~ 2 \ss (b, c ) \leq t < \infty .$ 
  			
  			Then $ \mathfrak{P}_{\theta} ( a, c, s ) =  	\dfrac{100 ~ \ss ( a, c ) }{ s } $ and $  \mathfrak{P}_{\theta} ( b, c, t ) =  	\dfrac{100 ~ \ss ( b, c ) }{ t } $. 
  			
  			Since $s + t \geq  2 \ss ( a, c ) + 2 \ss ( b, c ) \geq 2 \ss ( a, b ) $, hence $ \mathfrak{P}_{\theta} ( a, b, s + t ) =   	\dfrac{100 ~ \ss ( a, b ) }{ s + t  }  $ and
  			\begin{align*}
  				\mathfrak{P}_{\theta} ( a, b, s + t ) & =   	\dfrac{100 ~ \ss ( a, b ) }{ s + t  } 
  				\\
  				& \leq   100  \left(  \dfrac{ \ss ( a, c ) + \ss ( b, c ) }{ s + t  }    \right) 
  				\\
  				& < 100 \left( \dfrac{ \ss ( a, c ) }{ s }  +  \dfrac{  \ss ( c, b ) }{  t  } \right)  
  				=   2 \cdot \theta \bigg( \mathfrak{P}_{\theta} ( a, c, s ),  \mathfrak{P}_{\theta} ( b, c, t ) \bigg).  
  			\end{align*}
  			\item  Let $  2 \ss ( a, c ) \leq s < \infty, ~~  \ss (b, c ) < t \leq 2  \ss (b, c ) .$ 
  			Then, $ s + t > 2 \ss ( a, c ) + \ss (b, c ) > \ss ( a, c ) + \ss (b, c ) \geq \ss ( a, b ) $ and hence $ \mathfrak{P}_{\theta} ( a, b, s + t ) \leq 50  $. Again, $ \mathfrak{P}_{\theta} ( a, c, s ) =  	\dfrac{100 ~ \ss ( a, c ) }{ s } $ and $ \mathfrak{P}_{\theta} ( b, c, t ) =  	50 $  implies $  \mathfrak{P}_{\theta} ( a, b, s + t ) <  \mathfrak{P}_{\theta} ( a, c, s ) + \mathfrak{P}_{\theta} ( b, c, t ) = 2 \cdot \theta ( \mathfrak{P}_{\theta} ( a, c, s ),  \mathfrak{P}_{\theta} ( b, c, t ) ). $
  			\item Let $  2 \ss ( a, c ) \leq s < \infty, ~~  0 < t \leq   \ss(b, c ) .$ 
  			Then $  \mathfrak{P}_{\theta} ( a, c, s ) =  	\dfrac{100 ~ \ss ( a, c ) }{ s } $ and $ \mathfrak{P}_{\theta} ( b, c, t ) =  	25 $. If  $ s + t > 2 \ss ( a, c )  > 0 $ then $ \mathfrak{P}_{\theta} ( a, b, s + t ) \leq 50  $. Here we arrived at two cases. 
  			\begin{enumerate}[(a)]
  				\item Firstly, $ \mathfrak{P}_{\theta} ( a, b, s + t ) \leq 25 $. 
  				Then $  \mathfrak{P}_{\theta} ( a, b, s + t ) <  \mathfrak{P}_{\theta} ( a, c, s ) + \mathfrak{P}_{\theta} ( b, c, t ) = 2 \cdot \theta ( \mathfrak{P}_{\theta} ( a, c, s ),  \mathfrak{P}_{\theta} ( b, c, t ) ). $ 
  				\item Secondly, $ 25 < \mathfrak{P}_{\theta} ( a, b, s + t ) \leq 50 $. 
  				Then $$  \mathfrak{P}_{\theta} ( a, c, s ) +  \mathfrak{P}_{\theta} ( b, c, t ) = 	\dfrac{100 d ( a, c ) }{ s } + 	25  \begin{cases}
  					\geq 25 + 25 = 50 ~ ~~~~ \mbox{if} ~~  25 \leq \dfrac{100 d ( a, c ) }{ s } \leq 50 \\
  					> 25 = \frac{1}{2} \cdot 50  ~~~~~~~ \mbox{if} ~~~  0 \leq \dfrac{100 d ( a, c ) }{ s } \leq 25 
  				\end{cases}
  				.$$
  				Therefore, if $ ~ 25 \leq \dfrac{100 \ss ( a, c ) }{ s } \leq 50  $, then $ \mathfrak{P}_{\theta} ( a, b, s + t ) \leq 2 \cdot \theta ( \mathfrak{P}_{\theta} ( a, c, s ),  \mathfrak{P}_{\theta} ( b, c, t ) ) .$ Also  if  $  0 \leq \dfrac{100 d ( a, c ) }{ s } \leq 25 $ then
  				\begin{align*} 
  					\mathfrak{P}_{\theta} ( a, b, s + t )  = 50 < 2 \cdot ( \mathfrak{P}_{\theta} ( a, c, s ) + \mathfrak{P}_{\theta} ( b, c, t ) )  
  					= 4  \cdot \theta ( \mathfrak{P}_{\theta} ( a, c, s ),  \mathfrak{P}_{\theta} ( b, c, t ) ). 
  				\end{align*}
  			\end{enumerate}
  			The other cases can be verified similarly.   
  		\end{enumerate}
  		In particular, if we take $ f (\xi ) = \ln \xi, ~ \xi > 0 $, then 
  		$$  (\mathfrak{P}_{\theta} 2): \ln ( \mathfrak{P}_{\theta} ( x, \sigma, s+ t ) ) <  \ln \left(  \theta (	\mathfrak{P}_{\theta} ( x, z, s ) ,  	\mathfrak{P}_{\theta} ( z, \sigma,  t ) )   \right) + \ln 4   $$
  		holds  for all $ x, \sigma, z \in \mathcal{X} $ and  for all $ s, t > 0 $.  
  		
  		Therefore, $ ( \mathcal{X}, \mathfrak{P}_{\theta}, \theta ) $ is a generalized $ \theta $-paramatric metric space 
  		w.r.t   $ ( \ln x, \ln 4) \in \mathfrak{F} \times [ 0, \infty ) $.
  	\end{eg}

  	\begin{rem}
  		We know every continuous binary operation $`o$'   is a $ \mathcal{B} $-action (see \cite{8}). Let  $ \mathcal{B} $-action $ \theta ( \eth_1, \eth_2 ) = \dfrac{\eth_1 + \eth_2 }{ 2 }  $ in Example \ref{inc-dec p}. It is not hard to see  $`o$' is not a binary operation. Therefore, $ \mathcal{B} $-action of \cite{1} instead of  binary operation $`o$' of \cite{8}, one can obtain a larger class of  distance functions.
  	\end{rem}

  	
  	In the following result, we   construct a generalized parametric metric that induces from a generalized $ \theta $-paramatric metric which  establish the  connection of a generalized $ \theta $-parametric metric to generalized parametric metrics under  a given binary operation.

  	\begin{ppn}
  		Let $`o$' be a binary operation defined as in \cite{8} and $ ( \mathcal{X}, \mathfrak{P}_{\theta}, o ) $ be a generalized $ \theta $-parametric metric space with respect to a pair $ ( f, \alpha ) \in \mathfrak{F} \times [ 0, \infty ) $. Define a function   $ \mathcal{P} : \mathcal{X}^2 \times ( 0, \infty ) \to  [ 0, \infty ) $ by 
   $ \mathcal{P} ( x, \eta, t ) = \inf \bigg\{ 	\mathfrak{P}_{\theta} ( x, z, t_1 )  ~ o ~  	\mathfrak{P}_{\theta} ( z, \eta,  t_2 ) : z \in \mathcal{X}; ~ t_1 > 0, ~ t_2 > 0 ~ \& ~ t = t_1 + t_2 \bigg\}  $  
  		for all $ x, \eta \in \mathcal{X} $ and $ t > 0 $. If $ \mathfrak{P}_{\theta} $ satisfies  
  		$ ( \mathfrak{P}_{\theta} 3 ): ~ \mathfrak{P}_{\theta} ( x, \eta, t ) ~ \mbox{is  a  non-increasing  fuction  of} ~ t, ~ \mbox{for all} ~ x, \eta \in \mathcal{X} $ then $ \mathcal{P}  $ is a generalized parametric metric on $\mathcal{X} $ with respect to binary operation $`o.$' 
  		\begin{proof}
  			Let $ ( f, \alpha ) \in \mathfrak{F} \times [ 0, \infty ) $ satisfies $ ( \mathfrak{P}_{\theta} 2) $ with respect to which  $ \mathfrak{P}_{\theta} $. One can see   
  			\begin{enumerate}[(i)]
  				\item Since $ \mathfrak{P}_{\theta} ( x, x,  q ) = 0 $ for all $ q > 0 .$, By the definition of $ \mathcal{P} $,  $ \mathcal{P} ( x, x,  q ) = 0 $ for all $ q > 0 $. Let  $ x, \sigma \in \mathcal{X} $ with $ x \neq \sigma $ and suppose $ \mathcal{P} ( x, \sigma, t ) = 0 $ for all $ t > 0 $. Since $ x \neq \sigma $, then there exists $ t_0 > 0 $ such that $ \mathfrak{P}_{\theta} ( x, \sigma, t_0 ) > 0 $. Let $ \epsilon > 0 $. From the definition of $ P $, there exists $ z \in \mathcal{X} $ such that 
  				\begin{align*}
  					&	\theta \bigg( \mathfrak{P}_{\theta} ( x, z, t_1 ), \mathfrak{P}_{\theta} ( \sigma, z, t_2 ) \bigg) < \epsilon ~~~~ \mbox{where} ~ t_0 = t_1, t_2; ~ t_1>0, ~ t_2 > 0 \\
  					\implies & f \bigg( \theta ( \mathfrak{P}_{\theta} ( x, z, t_1 ), \mathfrak{P}_{\theta} ( \sigma, z, t_2 ) ) \bigg) \leq f (  \epsilon ) ~~~~ (\mbox{using} ~ (\mathfrak{F}_1))  \\
  					\implies & f \bigg( \mathfrak{P}_{\theta} ( x, \sigma, t_0 ) \bigg) \leq f (  \epsilon ) + \alpha ~~~ (\mbox{using} ~ (\mathfrak{P}_{\theta} 2)).   
  				\end{align*}
  			
  				Taking limit as  $ \epsilon \to 0^+ $,    we get 
  				\begin{align*}
  					&	\underset{ \epsilon \to 0^+}{lim} ~ ( f (  \epsilon ) + \alpha  ) = - \infty \\
  					or ~~ & f ( \mathfrak{P}_{\theta} ( x, \sigma, t_0 ) ) = - \infty \\
  					or ~~&  \mathfrak{P}_{\theta} ( x, \sigma, t_0 ) = 0 ~~~ (using ~ (\mathfrak{F}_2)).
  				\end{align*}
  				This is a contradiction to the assumption. Hence $ \mathcal{P} ( x, \sigma, t ) = 0 $ for all $ t > 0 $ implies $ x = \sigma $. Therefore $ (\mathcal{P}1) $ holds for $ \mathcal{P} $. 
  				\item Since $ \mathfrak{P}_{\theta} ( x, \eta, t ) = \mathfrak{P}_{\theta} (  \eta, x,  t ) $ for all $ t > 0 $ and for all $ x, \eta \in \mathcal{X} $, thus the function $ \mathcal{P} $  satisfies $ (\mathcal{P}2) $.   
  				\item Let $ x, y, z \in \mathcal{X} $ and   $ \epsilon > 0 $. By definition of $ \mathcal{P} $, there exist two chains $ x, u, y $ and  $ y, v, z $ such that
  				$$ \theta ( \mathfrak{P}_{\theta} (x, u, t_1 ), \mathfrak{P}_{\theta} ( u, y, t_2 ) ) < \mathcal{P} ( x, y, t ) + \epsilon ~ ~ or ~~ \mathfrak{P}_{\theta} (x, u, t_1 ) ~ o ~ \mathfrak{P}_{\theta} ( u, y, t_2 )  < \mathcal{P} ( x, y, t ) + \epsilon $$
  				and 
  				$$  \theta ( \mathfrak{P}_{\theta} (y, v, s_1 ), \mathfrak{P}_{\theta} ( v, z, s_2 ) ) < \mathcal{P} ( z, y, s ) + \epsilon ~~ or ~~ \mathfrak{P}_{\theta} (y, v, s_1 ) ~ o ~ \mathfrak{P}_{\theta} ( v, z, s_2 )  < \mathcal{P} ( z, y, s ) + \epsilon $$
  				where $ t = t_1 + t_2 $ and $ s = s_1 + s_2 $. Using the inequality $ (\mathfrak{P}_{\theta} 2) $ and the definition of $ \mathcal{P} $, we get 
  				\begin{align*}
  					\mathcal{P} ( x, z, s + t ) \leq & ~ \mathfrak{P}_{\theta} ( x, y, t ) ~ o ~ \mathfrak{P}_{\theta} ( z, y, s )  \\
  					\leq & ~ \{ \mathfrak{P}_{\theta} (x, u, t_1 ) ~ o ~ \mathfrak{P}_{\theta} ( u, y, t_2 ) \} ~ o ~ \{ \mathfrak{P}_{\theta} (y, v, s_1 ) ~ o ~ \mathfrak{P}_{\theta} ( v, z, s_2 ) \}  \\
  					< & ~ \left( \mathcal{P} ( x, y, t ) + \epsilon  \right) ~ o ~ \left(  P ( z, y, s ) + \epsilon \right).
  				\end{align*}
  				Allowing $ \epsilon \to 0^+ $,  
  				$ 	\mathcal{P} ( x, z, s + t ) \leq   \mathcal{P} ( x, y, t )   ~ o ~    \mathcal{P} ( z, y, s ). $ Thus $ \mathcal{P} $ satisfies $ ( \mathcal{P}3) $. 
  			\end{enumerate}	
  			Hence $ ( \mathcal{X}, \mathcal{P}, o ) $ is a generalized  parametric metric space. 
  		\end{proof}
  	\end{ppn}

  	\subsection{Basic properties of generalized $ \theta $-parametric metric space }
  	
  	This Section is divided into two subsections. In the first Section, we define a topology induced by  generalized $ \theta $-parametric metric, open \& closed ball and some other characterization.
  	The other section, consists of notion of convergent and Cauchy sequence and results related to their relations.  
  	
  	\subsection{Topological properties of  generalized $ \theta $-parametric metric space:}
  	
  	We start this part of topological properties of generalized $ \theta $-parametric metric space with the  definition of open ball and closed.
  	
  	\begin{dfn}
  		For each $ \eth > 0 $ and $ q > 0 $, we define open ball $ \mathscr{B} ( x, \eth, q ) $ and closed ball    $  \mathscr{B} [ x, \eth, q ] $ with center $ x \in \mathcal{X} $ and radius $ \eth > 0 $   as
  		$$  \mathscr{B} ( x, \eth, q )  = \bigg\{ y \in \mathcal{X} : \mathfrak{P}_{\theta} ( x, y, q ) < \eth \bigg\} ~~ \& ~~  \mathscr{B} [ x,\eth, q ] = \bigg\{ y \in \mathcal{X} : \mathfrak{P}_{\theta} ( x, y, q ) \leq \eth \bigg\}. $$
  	\end{dfn}

  	\begin{thm}
  		In a generalized $ \theta $-parametric metric space $ (\mathcal{X}, \mathfrak{P}_{\theta}, \theta ) $, the set 
  		$$ \tau_\theta = \{ \mathfrak{A} \subseteq \mathcal{X} : \text{for each} ~ a \in \mathfrak{A}  ~\text{and} ~ q > 0, ~ \text{there exists} ~ \eth > 0   ~ \mbox{such ~that} ~ \mathscr{B} ( a, \eth, q ) \subseteq \mathfrak{A} \} $$ 
  		is a topology on $\mathcal{X}$.
  	\end{thm}
  	\begin{proof}
  	 \begin{enumerate}[(i)]
  	 	\item Clearly, $ \emptyset, \mathcal{X} \in \tau_\theta $.
  	 	\item Let $ \mathscr{A}_1, \mathscr{A}_2 \in \tau_\theta $ and $ a \in \mathscr{A}_1 \cap \mathscr{A}_2 $. 
  	 	
  	 	So if we choose some $ t'   $ arbitrarily in $ ( 0, \infty ) $, then there exist  $ \eth_1, \eth_2 > 0 $ such that 
  	 	$ \mathscr{B} (a, \eth_i, t'  ) \subset \mathscr{A}_i $, for $  i =1, 2 $. 
  	 	
  	 	Let  $ \eth = \min \{ \eth_1, \eth_2 \} $. Then 
  	 	$$     x \in \mathscr{B} (a, \eth, t_i ) \implies P ( a, x, t' ) < \eth \leq \eth_i \implies x \in \mathscr{B} (a, \eth_i, t' ), ~~ for ~  i= 1, 2.   $$
  	 	Therefore, $   \mathscr{B} (a, \alpha, t' )  \subseteq  \mathscr{B} (a, \alpha_i, t' ) \subset \mathscr{A}_i, ~ i =1, 2 $, which implies    $  \mathscr{B} (a, \eth, t ) \subset \mathscr{A}_1 \cap \mathscr{A}_2 $. Hence   $ \mathscr{A}_1 \cap \mathscr{A}_2 \in \tau_\theta $. 
  	 	\item Let $  \mathscr{A}_i \in \tau_\theta, ~ i \in \Delta $ and $ a \in \cup_{ i \in \Delta } \mathscr{A}_i $. 
  	 	
  	 	Then  $ a \in \mathscr{A}_j $,  for some $ j \in \Delta $.  
  	 	Hence for some arbitrary $ t'   \in   ( 0, \infty ) $,  there exists $  \alpha > 0 $ such that 
  	 	$ B(a, \alpha, t) \subset \mathscr{A}_j \subset \cup_{ i \in \Delta } \mathscr{A}_i $.   
  	 	Therefore  	$ \cup_{ i \in \Delta } \mathscr{A}_i \in \tau_\theta $.
  	 \end{enumerate} 
  	\end{proof}

 \begin{rem}
  The members of $ \tau_\theta $ are called open sets in $ (\mathcal{X}, \mathfrak{P}_{\theta}, \theta ) $. 
 \end{rem}

  	\begin{thm}
  		Let $ ( \mathcal{X}, \mathfrak{P}_{\theta}, \theta ) $ be a generalized   $ \theta $-parametric metric  space. Then $ \tau_\theta $ is a Hausdorff topology on $\mathcal{X}$.
  	\end{thm}
  	\begin{proof}
  		Let $ ( f, \alpha ) \in \mathfrak{F} \times [ 0, \infty ) $ with respect to which  $ \mathfrak{P}_{\theta} $ satisfies $ ( \mathfrak{P}_{\theta} 2) $. Let $ x, \sigma \in \mathcal{X} $ be such that $ x \neq \sigma $. Then there exists $ t _0 > 0 $ such that $ \mathfrak{P}_{\theta} ( x, \sigma, t_0 ) \neq 0 $. We take $  a_\iota = \frac{\mathfrak{P}_{\theta} ( x, \sigma, t_0 )}{\iota }, ~ \forall \iota $. Then $ a_\iota \to 0 $ as $ \iota \to \infty $, by $( \mathfrak{F}_2) $: 
  		$ f ( a_\iota ) \to - \infty $ whenever $ \iota \to \infty $.
  		Since $ a_\iota > 0 $ for each $ \iota \in \mathbb{N} $, so using the axiom $ (\mathcal{B}3) $
  		for each $ k \in \mathbb{N} $, $ 0 < b_k < a_k $, there exists $ 0 < c_k < a_k $ such that $ \theta ( b_k, c_k ) = a_k $.
  		Now we consider two open sets  $ \mathfrak{U} = int ( \mathscr{B} ( x, b_k, \frac{t_0}{2}  )) $ and $  \mathfrak{V} = int( \mathscr{B} ( \sigma, c_k, \frac{t_0}{2} )  ) $  containing $ x $ and $ \sigma $ respectively.
  		Finally we claim  that $ \mathfrak{U} \cap \mathfrak{V} = \phi $.  
  		If possible suppose,  there exists $ z \in  \mathfrak{U} \cap \mathfrak{V} $. Then 
  		$$ \mathfrak{P}_{\theta} ( x, z, \frac{t_0}{2}   ) < b_k ~~~ \mbox{and} ~~~  \mathfrak{P}_{\theta} ( \sigma, z, \frac{t_0}{2}   ) < c_k. $$
  		Since $ \mathfrak{P}_{\theta} ( x, \sigma, t_0 ) > 0 $.  From the inequality  $ ( \mathfrak{P}_{\theta} 2) $, we have
  		$$ 	f (   \mathfrak{P}_{\theta}  ( x, \sigma, t_0 ) ) \leq   f \left( \theta (  \mathfrak{P}_{\theta} ( x, z, \frac{t_0}{2}   ),  \mathfrak{P}_{\theta} ( \sigma, z, \frac{t_0}{2}   )  ) \right) + \alpha  
  		<  f \left( \theta (   b_k, c_k ) \right) + \alpha  
  		=  f ( a_k ) + \alpha.  	$$ 
  		Taking limit as $ k \to \infty $, we get $ f ( a_k ) + \alpha \to - \infty ~~ or ~~ f (   \mathfrak{P}_{\theta}  ( x, \sigma, t_0 ) ) = - \infty  .$ Hence $  \mathfrak{P}_{\theta}  ( x, \sigma, t_0 ) = 0 $.  This contradicts our assumption.
  	\end{proof}

  	Next we present an example to see an open ball may not be an open set  in a generalized $ \theta $-parametric metric space.

  	\begin{eg}\label{ex via b-metric}
  		Let $ (\mathcal{X}, \ss ) $ be a b-metric space with constant coefficient $ K (\geq 1 ).$  Let us define   $ \mathfrak{P}_{\theta} : \mathcal{X}^2 \times ( 0, \infty ) \to  [ 0, \infty ) $  by $ \mathfrak{P}_{\theta} ( x, \sigma, q ) = \frac{ \ss ( x, \sigma ) }{q} ~~ \forall x, \sigma \in \mathcal{X} ~\& ~   q > 0.$
  		Then $ \mathfrak{P}_{\theta} $ satisfy  axiom $ ( \mathfrak{P}_{\theta} 1) $. Moreover, if for all $ \eth_1, \eth_2 \in [ 0, \infty ) $ we consider $ \mathcal{B} $-action be $ \theta ( \eth_1, \eth_2 ) =  { \eth_1 + \eth_2 }, $ we have 
  		\begin{align*}
  			& \mathfrak{P}_{\theta} ( x, \sigma, s + t )  \leq \frac{\ss ( x, \sigma ) }{ s + t }  
  			\leq   K \cdot \left( \frac{ \ss ( x, z ) + \ss ( z, \sigma ) }{ s + t } \right) 
  			<  K \cdot \left( \frac{ \ss ( x, z )   }{ s   }  + \frac{  \ss ( z, \sigma ) }{   t }  \right) \\
  			or ~~ & \ln ( \mathfrak{P}_{\theta} ( x, \sigma, s + t )) < \ln K + \ln (   \mathfrak{P}_{\theta} ( x, z, s ) +  \mathfrak{P}_{\theta} ( \sigma, z, t ) ). 
  		\end{align*} 
  		
  		Thus for $ f ( \eta ) = \ln \eta, ~ \eta > 0 $ and $ \alpha = \ln K ( \geq 0 ) $, $ P $ satisfy the axiom $ (\mathfrak{P}_{\theta} 2) .$  Therefore $ \mathfrak{P}_{\theta} $ is a generalized $ \theta $-parametric metric   on $\mathcal{X}$. 
  		
  		Next let $ \mathcal{X} = \{ 0, 1, \frac{1}{2}, \cdots \} .$ Let us consider b-metric space $ ( \mathcal{X}, \ss ) $ of \cite[Example 3.9]{14} with    
  		$$ \ss ( \xi, \eta ) = \begin{cases}
  			0 \quad \hskip 20 pt ~ \mbox{if} ~  \xi = \eta \\
  			1 \quad \hskip 20 pt ~ \mbox{if} ~  \xi \neq \eta \in \{ 0, 1 \} \\
  			| \xi - \eta | \quad ~\mbox{if} ~  \xi \neq \eta \in \{ 0 \} \cup \{ \frac{1}{2m} : m = 1, 2, \cdots  \} \\
  			4 \quad  \hskip 20 pt ~~ \mbox{otherwise}
  		\end{cases}  $$
  		and $ K = \frac{8}{3} .$ Then the open ball $ \mathscr{B} ( 1, 2, 1 ) $ is the set \begin{align*}   
  			\mathscr{B} ( 1, 2, 1 )  = \bigg\{ y \in \mathcal{X} : \mathfrak{P}_{\theta} ( y, 1, 1 ) < 2 \bigg\}  
  			= \bigg\{ y \in \mathcal{X} : \ss ( y, 1 ) < 2 \bigg\}  
  			= \{ 0, 1 \} 
  		\end{align*}
  		in $  ( \mathcal{X}, \mathfrak{P}_{\theta}, + ) .$  
  		Thus $ 0 \in  \mathscr{B} ( 1, 2, 1 ) .$ Let $ \eth_0 > 0 $ be fixed.
  		We  claim that there  does not  exist any  $ t_0 > 0 $ such that $ \mathscr{B}( 0, \eth_0, t_0 ) \subset  \mathscr{B} ( 1, 2, 1 ) $. 	
  		For the fixed $ \eth_0 > 0 $ and for any $ t_0 > 0 $, by Archimedean property we can find $ N \in \mathbb{N} $ such that $	0 < \frac{1}{(2 N) t_0 } < \eth_0 $ or $  \dfrac{\ss ( 0, \frac{1}{2 N } )}{t_0}  < \eth_0. $ Clearly 
  		$  \frac{1}{2 N } \in \mathscr{B} ( 0, \eth_0, t_0 ) $.  However $ \mathfrak{P}_{\theta} ( 1,  \frac{1}{2 N }, 1   ) = \frac{\ss ( 1, \frac{1}{2 N } ) }{1} = 4 \nleq 2 $. Hence $ \frac{1}{2 N } \notin \mathscr{B} ( 1, 2, 1 ) .$ Consequently $ \mathscr{B} ( 0, \eth_0, t_0 ) \nsubseteq  \mathscr{B} ( 1, 2, 1 ) .$   Therefore $ \mathscr{B} ( 1, 2, 1 ) $     is not an  open set in $  (\mathcal{X}, \mathfrak{P}_{\theta}, + ) .$   
  	\end{eg}

  	Next, we search a sufficient condition such that an open ball is a possible to be an open set in generalized $ \theta $-parametric metric space.
  	
  	 To fill our admire, let $ \eth > 0 $ be fixed. For some $ t > 0 $ and $ x \in \mathcal{X} $, we consider open ball $ \mathscr{B} ( x, \eth, t ) $ is in $ ( \mathcal{X}, \mathfrak{P}_{\theta}, \theta ) $.  Suppose  $   \mathfrak{P}_{\theta}  $ satisfy  axiom  $  ( \mathfrak{P}_{\theta} 2) $ whenever  $ ( f, \alpha ) \in \mathfrak{F} \times [ 0, \infty  ) .$ Then $ f ( \eth ) - \alpha \in \mathbb{R} .$  Further, by  $ ( \mathfrak{F}_2) $, there exists $ \delta > 0 $ such that 
  	\begin{equation} \label{open ball open set 1}
  		0 < s < \delta ~ \implies ~ f (s ) < f ( \eth )  - \alpha. 
  	\end{equation}
  	Clearly $ \delta $ depends on $ \eth $ and $ \alpha .$ We denote this expression by $h ( \eth, \alpha ) $. We are ready to find a sufficient condition of open ball to be an open set in generalized $\theta$-parametric metric space with the help of $h ( \eth, \alpha ) $ in the following theorem.

  	\begin{thm}
  		Let $ ( \mathcal{X}, \mathfrak{P}_{\theta}, \theta ) $ be  a generalized $ \theta $-parametric metric space with respect to a pair $ ( f, \alpha ) \in \mathfrak{F} \times [ 0, \infty  ) .$  Consider for all $ a, b \in \mathcal{X}  $, $ \mathfrak{P}_{\theta} ( a, b, u ) $ be a continuous function of $ u > 0 .$ If $ \eth < h ( \eth, \alpha ) ,$ then the open ball  $ \mathscr{B} ( x, \eth, t ) $ is  an   open set. 
  	\end{thm}
  	\begin{proof}
  		Let $ y \in  \mathscr{B} ( x, \eth, t ) .$ Then $ \mathfrak{P}_{\theta} ( x, y, t) < \eth $.  Since $  \mathfrak{P}_{\theta} ( x, y, u ) $  is  a continuous function of $ u >  0 $, there exists $t _0 \in ( 0, t ) $ so that $  \mathfrak{P}_{\theta} ( x, y, t_0 ) < \eth.$ Suppose $ \mathfrak{P}_{\theta} ( x, y, t_0 ) = \eth_0 (< \eth ) .$ Using condition $ ( \mathscr{B} 3) $, we can find $ \eth_1 \in ( 0, \eth ) $ such that $ \theta ( \eth_0, \eth_1 ) = \eth $.  We consider the open ball $ \mathscr{B} ( y, \eth_1, t_1 ) $ whenever $ t_1 = t - t_0 .$  Let $ z \in  \mathscr{B} ( y, \eth_1, t_1 ) .$ In this occasion $ \mathfrak{P}_{\theta} ( z, y, t_1) < \eth_1 ,$  and we have
  		\begin{align*}  
  			\theta ( \mathfrak{P}_{\theta} ( x, y, t_0 ),  \mathfrak{P}_{\theta} ( z, y, t_1) ) &< \theta ( \eth_0, \eth_1 ) \\
  			& = \eth < h ( \eth, \alpha ). 
  		\end{align*} 
  		This gives us $f	\left(  \theta ( \mathfrak{P}_{\theta} ( x, y, t_0 ),  \mathfrak{P}_{\theta} ( z, y, t_1) )  \right)   < f ( \eth )  - \alpha$ whenever we use (\ref{open ball open set 1}). Further, by ($\mathfrak{P}_{\theta} 2)$,	$f ( \mathfrak{P}_{\theta} ( x, z, t ) )  < f ( \eth ).$ Using the proprty   $(\mathfrak{F}_1)$ we have  $\mathfrak{P}_{\theta} ( x, z, t )   <   \eth.$   Hence $ z \in  \mathscr{B} ( x, \eth, t ) .$ Consequently, $  \mathscr{B} ( y, \eth_1, t_1 )  \subset \mathscr{B} ( x, \eth, t ) .$ Therefore $  \mathscr{B} ( x, \eth, t ) $ is an open set in   
  	  $ ( \mathcal{X}, \mathfrak{P}_{\theta}, \theta ) $. 
  	\end{proof}

 \begin{rem}
  Although we prove that an open ball is an open set  by providing a sufficiet condition on   an upper bound for the radius, determination of a necessary and sufficient condition under which an open ball in a generalized $\theta$-parametric metric space becomes an open set still remains an interesting problem for future investigation.
 \end{rem}

  	

  	\subsection{Convergence of sequence in   generalized $ \theta $-parametric metric space:}
  	
  	This section explores the fundamental properties of the generalized $\theta$-parametric metric space, focusing on Cauchy and convergent sequences and their interrelations. We begin with the   definition of convergent sequence.

  	\begin{dfn} \label{conv seq dfn}
  		Let $ \{ s_\iota \} $ be a sequence in a generalized $ \theta $-parametric metric space $ ( \mathcal{X}, \mathfrak{P}_{\theta}, \theta ) $.
  		\begin{enumerate}[(i)]
  			\item  $ \{ s_\iota \} $ is converge to   $ s \in \mathcal{X} $, if for a given  $ \epsilon > 0 $,   for each $ p > 0 $, there exists a natural number $ N ( p )   $ such that $  \mathfrak{P}_{\theta}  ( s_\iota, s, p ) < \epsilon $ for all $ \iota \geq N( p ). $
  			\item   $ \{ s_\iota \} $ is a Cauchy sequence   if for a given  $ \epsilon > 0 $,   for each $ p > 0 $, there exists a natural number $ N (p )   $  such that $  \mathfrak{P}_{\theta}  (  s_{\iota_1}, s_{\iota_2}, p ) < \epsilon $ for all $ {\iota_1}, {\iota_2} \geq N( p ). $
  		\end{enumerate}
  	\end{dfn}

  	The following results are immediate consequence of Definition \ref{conv seq dfn}.

  	\begin{ppn}
  		Let $ \{ s_n \} $ be a sequence in a generalized $ \theta $-parametric metric space $ ( \mathcal{X}, \mathfrak{P}_{\theta}, \theta ) $. Then 
  		\begin{enumerate}[(i)]
  			\item $ \{ s_n \} $   converges to $ s $ iff $ \underset{\iota \to \infty }{lim }  ~ \mathfrak{P}_{\theta}  ( s_\iota, s, p ) = 0 ~$ for all $ p > 0 $. 
  			\item $ \{ s_n \} $ is a Cauchy sequence iff $ \underset{\iota_1, \iota_2 \to \infty }{lim }  ~ \mathfrak{P}_{\theta}  ( s_{\iota_1}, s_{\iota_2}, p ) = 0 ~$ for all $ p > 0 $. 
  		\end{enumerate}
  	\end{ppn}
  	\begin{proof}
  		The proof follows from the Definition \ref{conv seq dfn} and properties $ (\mathfrak{F}_1)$-$(\mathfrak{F}_2) $.
  	\end{proof}

  	\begin{ppn} \label{conv-cauchy relation}
  		Let $ \{ s_\iota \} $ be a convergent sequence in a generalized $ \theta $-paramatric metric space $ ( \mathcal{X}, \mathfrak{P}_{\theta}, \theta ) $.
  		\begin{enumerate}[(i)]
  			\item Converging point of $ \{ s_\iota \} $ is unique  in $ \mathcal{X} $. 
  			\item $ \{ s_\iota \} $ is a Cauchy sequence in $ \mathcal{X} $. 
  		\end{enumerate}
  	\end{ppn}
  	\begin{proof}
  		Let $ ( f, \alpha ) \in \mathfrak{F} \times [ 0, \infty ) $ be the pair with respect to which $ \mathfrak{P}_{\theta} $ satisfies $ ( \mathfrak{P}_{\theta} 2) $.
  		\begin{enumerate} [(i)]
  			\item If possible suppose that $ \{ s_\iota \} $ converges to $ x $ and $ y $ in $ \mathcal{X} $ where $x \neq y.$  Then 
  			\begin{equation}\label{uniquelimit 1}
  				\lim\limits_{\iota \to \infty } \mathfrak{P}_{\theta}  ( s_\iota, x, p ) = 0 =  \lim\limits_{\iota \to \infty }  ~ \mathfrak{P}_{\theta}  ( s_\iota, y, p ) = 0 ~~ \text{for all} ~ p > 0. 
  			\end{equation} 
  			Since $ x \neq y $, thus there exists $ p_0 > 0 $ such that 
  			\begin{equation}\label{uniquelimit 2}
  				\mathfrak{P}_{\theta} ( x, y, p_0 ) > 0.
  			\end{equation}
  		
  			Using $ (\mathfrak{P}_{\theta} 2) $,  we have 
  			$ f ( \mathfrak{P}_{\theta} ( x, y, p_0 ) ) \leq f \left( \theta ( \mathfrak{P}_{\theta} (x, s_\iota, \frac{p_0}{2} ), \mathfrak{P}_{\theta} ( s_\iota, y, \frac{p_0}{2} )   ) \right) + \alpha  ~~ \text{for all} ~ \iota \in \mathbb{N}. $
  			Using (\ref{uniquelimit 1}), we get 
  			\begin{align*}
  				& \underset{\iota \to \infty }{lim }  ~ \theta \left( \mathfrak{P}_{\theta} \left( x, s_\iota, \frac{p_0}{2} \right), ~ \mathfrak{P}_{\theta}   \left( s_\iota, y, \frac{p_0}{2} \right)   \right) = 0 \\
  				\implies & \underset{\iota \to \infty }{lim }  ~ f \left( \theta \left( \mathfrak{P}_{\theta} (x, s_\iota, \frac{p_0}{2}  ), \mathfrak{P}_{\theta} ( s_\iota, y, \frac{p_0}{2} )   \right) \right) + \alpha = - \infty  \\
  				\implies & \underset{\iota \to \infty }{lim }  ~ f ( \mathfrak{P}_{\theta} ( x, y, p_0 ) ) = - \infty  \\
  				\implies & \mathfrak{P}_{\theta} ( x, y,p_0 ) = 0 ~~~ ( \text{using} ~ (\mathscr{F}2))
  			\end{align*}
  		
  			This is a contradiction to the relation (\ref{uniquelimit 2}). 
  			\item Suppose the sequence $ \{ s_\iota \} $ convergs to   $ x  \in  \mathcal{X} $. Then 
  			$	\underset{\iota \to \infty }{lim }  ~ \mathfrak{P}_{\theta}  ( s_\iota, x, p ) = 0 ~~ \mbox{for~ all} ~ p > 0. $ Let $ \epsilon > 0 $. Using $ ( \mathfrak{F}_2) $, there exists $ \delta > 0 $ so that
  			\begin{equation}\label{conv cauchy 1}
  				0 < t < \delta \implies f( t ) < f ( \epsilon ) - \alpha. 
  			\end{equation}
  		
  			Next for $ \delta > 0 ,$ and for each $  p > 0 $, there exists $ N  ( p )  \in  \mathbb{N} $ we have 
  			$  \mathfrak{P}_{\theta} ( s_\iota, x, \frac{p}{2} ) <  \frac{ \delta}{2} ~~~ \text{for all} ~ n\iota \geq N  ( p ). $
  			Therefore for each $ p > 0 $, 
  			\begin{align*}
  				&  \theta \left( \mathfrak{P}_{\theta} ( s_{\iota_1}, x, \frac{p}{2} ), \mathfrak{P}_{\theta} ( x, s_{\iota_2}, \frac{p}{2} )   \right) < \theta ( \frac{ \delta}{2}, \frac{ \delta}{2} )  ~~ ~ \text{for all} ~ \iota_1, \iota_2 \geq N  ( p ) \\
  				\implies  & \theta \left( \mathfrak{P}_{\theta} ( s_{\iota_1}, x, \frac{p}{2} ), \mathfrak{P}_{\theta} ( x, s_{\iota_2}, \frac{p}{2} )   \right) \leq \frac{ \delta}{2} <   \delta ~~ ~ \text{for all} ~ \iota_1, \iota_2 \geq N  ( p ) \\
  				\implies  &  f \left( \theta \left( \mathfrak{P}_{\theta} (s_{\iota_1}, x, \frac{p}{2} ), \mathfrak{P}_{\theta} ( x, s_{\iota_2}, \frac{p}{2} )   \right) \right) < f ( \epsilon ) - \alpha ~~ ~ \text{for all} ~ \iota_1, \iota_2 \geq N  ( p ) \\
  				\implies  &  f ( \mathfrak{P}_{\theta} ( s_{\iota_1}, s_{\iota_2}, p ) ) < f ( \epsilon ) ~~ ~ \text{for all} ~ \iota_1, \iota_2 \geq N  ( p ) \\
  				\implies  &  \mathfrak{P}_{\theta} ( s_{\iota_1}, s_{\iota_2}, p ) < \epsilon ~~ ~ \text{for all} ~ \iota_1, \iota_2 \geq N  ( p ). 
  			\end{align*}
  		
  			This proves that $\{ s_\iota \} $ is a Cauchy sequence in $\mathcal{X}$.  
  		\end{enumerate}
  	\end{proof}

\begin{dfn}
	A generalized $ \theta $-paramatric metric space $ ( \mathcal{X}, \mathfrak{P}_{\theta}, \theta ) $   is said to be complete if every Cauchy sequence in $\mathcal{X}$ converges to some element in it. 
\end{dfn}

  The  Example \ref{ex via b-metric} highlights  fundamental deviation from classical metric space behavior. 
 This raises further questions about other fundamental properties induced by the generalized $\theta$-parametric metric.
 To explore this further, we turn our attention to the notion of sequential continuity.

 \begin{dfn}
 	 Let  $ ( \mathcal{X}, \mathfrak{P}_{\theta}, \theta ) $ be a generalized $ \theta $-paramatric metric space.  
 	The generalized $ \theta $-paramatric metric function $ \mathfrak{P}_{\theta}  $ is said to be sequentially continuous in its variables, if for any two sequence $ \{ \mu_\iota \} $ \& $ \{ \eta_\iota \} $    converging to $ \mu $ \& $ \eta $  respectively,  $ \underset{\iota \to \infty }{lim }  ~ \mathfrak{P}_{\theta}  (  \mu_\iota, \eta_\iota, p ) =   \mathfrak{P}_{\theta}  ( \mu, \eta, p )   $ holds for all $ p > 0 $. 
 \end{dfn}
 
 \begin{rem}
  In classical settings, the metric function	is typically sequentially continuous in both the variable. However, this property does not hold in general for the generalized $\theta$-parametric metric. We   provide a counterexample to illustrate this limitation. 
 	\begin{eg}
 	 Consider the generalized $\theta$-parametric metric space  of Example \ref{ex via b-metric}.  
 		Then observe that for the sequence $ \{ s_n \} = \{ \frac{1}{2n } \}$ in $ \mathcal{X} $, $ \underset{\iota \to \infty }{lim }  ~ \mathfrak{P}_{\theta}  ( s_\iota, 0, p ) = 0,~ \forall p > 0$ i.e. $ \{ \frac{1}{2n }  \} $ converges to $ 0 $. But, for all $ p > 0 $, 
 		  $ \underset{\iota \to \infty }{lim }  ~ \mathfrak{P}_{\theta}  ( s_\iota, 1, p ) =  \frac{4}{p} $ ~and~ $  \mathfrak{P}_{\theta}  ( 0, 1, p ) = \frac{1}{p  }  $, implies that  	$ \underset{\iota \to \infty }{lim }  ~ \mathfrak{P}_{\theta}  ( s_\iota, 1, p )  \neq \mathfrak{P}_{\theta}  ( 0, 1, p ), ~ \forall p > 0 $. 
 	\end{eg}
 \end{rem}

 To further examine the topological structure induced by  $\mathfrak{P}_\theta$, we now introduce the notion of a closed ball in this setting. While in classical metric spaces closed balls are always closed sets, this is not generally the case in our framework. The subsequent example demonstrates that a closed ball in a generalized $\theta$-parametric metric space need not be a closed set.

  	\begin{dfn}
  	  A subset $ F $ of a generalized $ \theta $-paramatric metric space  	 $ ( \mathcal{X}, \mathfrak{P}_{\theta}, \theta ) $   is said to be a closed set if there is no sequence in $F$ which converges in $ \mathcal{X} \setminus F$.
  	\end{dfn}

\begin{rem}
$ \overline{F} $ denotes the closure of   $ F  \subset \mathcal{X} $   with respect to the topology $ \tau_\theta $. 
\end{rem}

\begin{eg}
	\label{closed ball not closed set}
    Consider the b-metric space $ ( \mathcal{X}, \ss ) $ \cite[Example 3.10]{14} where $ \mathcal{X} =  \{ 0, 1, \frac{1}{2}, \cdots \} $ and $$ \ss ( \xi, \eta ) = \begin{cases}
	0 \quad \hskip 20 pt ~ \mbox{if} ~  \xi = \eta \\
	1 \quad \hskip 20 pt ~ \mbox{if} ~  \xi \neq \eta \in \{ 0, 1 \} \\
	| \xi - \eta | \quad ~\mbox{if} ~  \xi \neq \eta \in \{ 0 \} \cup \{ \frac{1}{2m} : m = 1, 2, \cdots  \} \\
	\frac{1}{4} \quad  \hskip 20 pt ~~ \mbox{otherwise}
\end{cases}  $$
with constant coefficiet $ K = 4  $.

 Then by   Example \ref{ex via b-metric}, we can say that $ \mathfrak{P}_{\theta} ( x, \sigma, q ) = \frac{ \ss ( x, \sigma ) }{q} ~~ \forall x, \sigma \in \mathcal{X} ~\& ~   q > 0 $ is a generalized $ \theta $-parametric metric   on $\mathcal{X}$ with respect to  $ f ( \eta ) = \ln \eta, ~ \eta > 0 $ and $ \alpha = \ln K ( \geq 0 ) $ under the $ \mathcal{B} $-action: $ \theta ( \eth_1, \eth_2 ) =  { \eth_1 + \eth_2 }  $. 
 
  Next obseve that the closed ball $ \mathscr{B} [ 1, \frac{1}{2}, 1 ] $ is the set \begin{align*}   
		\mathscr{B} [ 1, \frac{1}{2}, 1 ]   = \bigg\{ y \in \mathcal{X} : \mathfrak{P}_{\theta} ( y, 1, 1 ) \leq  \frac{1}{2} \bigg\}  
		= \bigg\{ y \in \mathcal{X} : \ss ( y, 1 )  \leq  \frac{1}{2} \bigg\}  
		= \mathcal{X} \setminus \{  0 \}. 
	\end{align*}
	 
	 Now consider the sequence $ \{ s_n \} = \{ \frac{1}{2n } \}$ in $ \mathcal{X} $.  Since $ \underset{\iota \to \infty }{lim }  ~ \mathfrak{P}_{\theta}  ( s_\iota, 0, p ) = \underset{\iota \to \infty }{lim } ~ \frac{1}{2 \iota  } = 0 ~$ for all $ p > 0 $, so $ \{ s_n \} $ converges to $ 0 $. But $ 0 \notin \mathscr{B} [ 1, \frac{1}{2}, 1 ]  $.   Therefore $\mathscr{B} [ 1, \frac{1}{2}, 1 ] $     is not a closed set in $  (\mathcal{X}, \mathfrak{P}_{\theta}, + ) .$   
 \end{eg}

 \begin{ppn} \label{closed ball closed set proof}
 		  Let $  \mathscr{B} [ a, \eth, t ] $ for some $a \in \mathcal{X} $ and $ \eth, t > 0 $, be  a closed ball  in a generalized $ \theta $-paramatric metric space   $ ( \mathcal{X}, \mathfrak{P}_{\theta}, \theta ) $. Suppose that for every sequence $\{ s_\iota \} \subset \mathcal{X} $, we have
 		$$ \underset{\iota \to \infty }{ \lim }  ~ \mathfrak{P}_{\theta}  ( s_\iota,  s, p ) = 0, ~ s \in \mathcal{X}, ~ \forall p > 0 \implies    \mathfrak{P}_{\theta}  ( s,  y, p ) \leq \underset{\iota \to \infty }{\lim \sup  } ~ ~ \mathfrak{P}_{\theta}  ( s_\iota,  y, p ), ~ y \in  \mathcal{X}, ~ \forall p > 0. $$
 		 
 Then  $  \mathscr{B} [ a, \eth, t ] $ is a closed set in $ ( \mathcal{X}, \mathfrak{P}_{\theta}, \theta ) $.
 \end{ppn}
\begin{proof}
	 Let $ \{ s_\iota \} \subset  \mathscr{B} [ a, \eth, t ] $ be a sequence converging to some $ s \in \mathcal{X} $. Then  
	 $$ \underset{\iota \to \infty }{ \lim }  ~ \mathfrak{P}_{\theta}  ( s_\iota,  s, p ) = 0, ~ \forall p > 0. $$   
	    
	    Again by the definition of $  \mathscr{B} [ a, \eth, t ] $, we have
	 $$  \mathfrak{P}_{\theta} ( s_\iota, a, t ) \leq \eth, ~ \forall \iota \in \mathbb{N}. $$ 
	 
	 Passing  supremum limit as $ \iota \to \infty $ on the above inequality, we get 
	 \begin{align*}
	 	 & \underset{\iota \to \infty }{\lim \sup  } ~ ~ \mathfrak{P}_{\theta}  ( s_\iota,  a, p ) \leq \eth \\
	 	 \implies & \mathfrak{P}_{\theta}  ( s,  a, p ) \leq  \eth \\
	 	 \implies & s \in \mathscr{B} [ a, \eth, t ]. 
	 \end{align*} 
 
  Therefore, $ \mathscr{B} [ a, \eth, t ] $ is a closed set.
\end{proof}

 \begin{rem}
     Proposition \ref{closed ball closed set proof} provides only a sufficient condition ensuring  that a  closed ball $ 	\mathscr{B} [ a, r, t ] $ is a closed set in a generalized $ \theta $-paramatric metric space. An interesting problem consists to find a necessary  and sufficient condition under which every closed ball becomes a closed set.
 \end{rem}

 \begin{ppn} \label{ppn for compact}
  Let $(\mathcal{X}, \mathfrak{P}_\theta, \theta)$ be a generalized $\theta$-parametric metric space   such that $ \mathfrak{P}_\theta ( a, b, u ) $ is a contiuous function of $ u > 0 $, for all $ a, b \in \mathcal{X}$.      Then   for any subset  $ F $ of $ \mathcal{X}$,
 $$ x \in \overline{F}, ~ r > 0 \implies \mathscr{B} ( x, r, t ) \cap F \neq \emptyset, \quad \text{for any} ~ t> 0. $$ 
  \end{ppn}
\begin{proof}
	Let $ ( f, \alpha ) \in \mathfrak{F} \times [ 0, \infty ) $ be the pair with respect to which $ \mathfrak{P}_{\theta} $ satisfies $ ( \mathfrak{P}_{\theta} 2) $. Let us define
	$$ \mathcal{O} = \{ x \in X : \text{for any } r > 0 ~ \& ~ t>0, ~ \exists a \in F ~\text{such that} ~ \mathfrak{P}_{\theta} ( x, a, t ) < r \}. $$
	
	Clearly, $ 	F \subset \mathcal{O} $. We claim that  $ \mathcal{O}  $ is a closed set in $ \mathcal{X}$. 
	
	Let $ \{ \xi_n \} $ be a sequence in  $ \mathcal{O} $ such that $ \{\xi_n\} $ converges to some $ x \in   \mathcal{X}$.  
	
	Choose $ r > 0 $ and $ t' > 0 $ arbitrarily. Then by $ (\mathcal{F}2) $, there exists $ \delta_r > 0 $ such that 
	\begin{equation}\label{12.4-1}
		0 < u <  \delta_r \implies f ( u ) < f ( r ) - \alpha. 
	\end{equation}

Since $ x_n \to x $, so for $ \mu ( \delta_r )$ ( a positive function of $ \delta_r > 0 $), there exists $ N \in \mathbb{N} $ such that 
\begin{equation}\label{12.4-2}
	 \mathfrak{P}_{\theta} ( \xi_n, x, s ) <  \mu ( \delta_r ) \quad \quad \forall n \geq N ~ \& ~ \forall s > 0.  
\end{equation} 

On the other hand, for $ \mu ( \delta_r ) > 0 $ and  $ \delta_r > 0 $, $ (\mathcal{B}3) $ ensures the existence of a $ \eta ( \delta_r ) $ ( a positive function of $ \delta_r > 0 $) such that 
$$ \theta ( \mu ( \delta_r ),  \eta ( \delta_r ) ) = \delta_r. $$

Again $ x_N \in \mathcal{O} $ implies there exists $a \in F$ such that
$$  \mathfrak{P}_{\theta} ( \xi_N, a, t' ) <    \eta ( \delta_r ). $$

Since $  \mathfrak{P}_{\theta} ( \xi_N, a, u ) $ is a contiuous function of $ u > 0 $, so there exists $ t'' \in ( 0, t' ) $ such that 
$$  \mathfrak{P}_{\theta} ( \xi_N, a, t'' ) <    \eta ( \delta_r ). $$

From (\ref{12.4-2}), in particular, we can write 
$$  \mathfrak{P}_{\theta} ( \xi_N, x, t' - t'' ) <  \mu ( \delta_r ). $$

Now $  \mathfrak{P}_{\theta} ( \xi_N, a, t' ) <    \eta ( \delta_r ) $ and  $  \mathfrak{P}_{\theta} ( \xi_N, x, t' - t'' ) <  \mu ( \delta_r ) $ together implies 
\begin{align} \label{12.4-3}
	 & \theta \bigg( \mathfrak{P}_{\theta} ( \xi_N, a, t' ) , \mathfrak{P}_{\theta} ( \xi_N, x, t' - t'' ) \bigg) < \theta ( \eta ( \delta_r ), \mu ( \delta_r ) ) = \delta_r \nonumber \\
	 \implies & f \bigg( \theta \bigg( \mathfrak{P}_{\theta} ( \xi_N, a, t' ) , \mathfrak{P}_{\theta} ( \xi_N, x, t' - t'' ) \bigg)  \bigg) < f ( r ) - \alpha \quad \quad (using ~ (\ref{12.4-1})).  
\end{align}

If $ \mathfrak{P}_{\theta} ( x, a, t' ) > 0 $, then by $ (\mathfrak{P}_{\theta}2 ) $, we get
\begin{align*}
	& f \bigg(   \mathfrak{P}_{\theta} ( x, a, t' )  \bigg) \leq f \bigg( \theta \bigg( \mathfrak{P}_{\theta} ( \xi_N, a, t' ) , \mathfrak{P}_{\theta} ( \xi_N, x, t' - t'' ) \bigg)  \bigg) + \alpha < f(r ) \quad (using ~ (\ref{12.4-3}))   \\
	 \implies & \mathfrak{P}_{\theta} ( x, a, t' )  < r \quad \quad (using ~ (\mathscr{F}1)) \\  
	 \implies & x \in \mathcal{O}. 
\end{align*}
 
 This proves that $ \mathcal{O} $ is a closed set in $ ( \mathcal{X}, \mathfrak{P}_{\theta}, \theta ) $, which contains $F$. Thus $ \overline{F} \subseteq \mathcal{O} $. This yields the conclusion.
\end{proof}

Next, we discuss the compactness on generalized $ \theta $-parametric metric space.

\begin{dfn}
  A subset $ F $ of a generalized $ \theta $-paramatric metric space 	$ ( \mathcal{X}, \mathfrak{P}_{\theta}, \theta ) $   is said to be a compact set if every open cover of $F$ has a finite subcover. 
\end{dfn}

\begin{lma} \label{lma for compact}
 	Let $(\mathcal{X}, \mathfrak{P}_\theta, \theta)$ be a generalized $\theta$-parametric metric space 	and   $ F \subseteq \mathcal{X}$ be compact. Then any decreasing sequence of nonempty closed subsets of $F$ has a nonempty intersection i.e. if $\{C_n\}_{n \in \mathbb{N}}$ is a sequence of closed subsets of $F $ such that
$	C_1 \supseteq C_2 \supseteq C_3 \supseteq \cdots \quad \text{and} \quad C_n \neq \emptyset \text{ for all } n $,  then
 $ 	\bigcap_{n=1}^\infty C_n \neq \emptyset. $
\end{lma}
\begin{proof} 
	The proof is similar to the proof of classical metric space. 
\end{proof}

In continuation, we present the following definition and proposition that   naturally comes after   the compactness.

\begin{dfn}
	  A subset $ F $ of a generalized $ \theta $-paramatric metric space 
	$ ( \mathcal{X}, \mathfrak{P}_{\theta}, \theta ) $   is said to be a compact set if   for any sequence sequence $ \{ s_\iota \} $ in $F  $,  there exist a subsequence $ \{ s_{\iota_k} \} $ of $ \{ s_\iota \} $ which converges to some element in $   F$.
\end{dfn}

  \begin{ppn}
   Let $ ( \mathcal{X}, \mathfrak{P}_{\theta}, \theta ) $  be a  generalized $ \theta $-paramatric metric space	such that $ \mathfrak{P}_\theta ( a, b, u ) $ is a contiuous function of $ u > 0 $, for all $ a, b \in \mathcal{X}$.   
	Then  for some subset  $ F $ of $  \mathcal{X} $,   $ F $ is compact iff   $ F $ is sequentially compact.
  \end{ppn}	
	\begin{proof}
		 First assume that $ F $ is a compact set in  $ \mathcal{X} $. 
		
		Let $ \{s_n\} $ be a sequence in $F$. Next we construct a sequence $ \{ B_\iota \} $ defined by
		$$ B_\iota  = \{ s_m : m \geq \iota  \}, \quad \iota  \in \mathbb{N}. $$
		
		Then clearly $ B_{n+1} \subset B_{n } $ and hence $  \{ \overline{  B_n } \} $ is a decreasing sequence of non-empty closed subsets of $F$. Therefore by Lemma \ref{lma for compact}, there exists $ x \in \underset{  \iota \in \mathbb{N}}{\cap} \overline{ B_\iota }  $ for certain $x\in F $.
		
		Let $ \epsilon $ be an arbitrary positive real number. 
		
		Since $ x \in   \overline{ B_1 }  $, so by Proposition \ref{ppn for compact},    for  any $ t > 0 $,
		\begin{align*}
			& \mathscr{B} ( x, \epsilon, t ) \cap F \neq \emptyset\\
			\implies & \exists n_1 \geq 1 ~ \& ~ s_{n_1} \in F ~ \text{such  that} ~ s_{n_1} \in  \mathscr{B} ( x, \epsilon, t ) \cap F \\
			\implies & \mathfrak{P}_\theta ( s_{n_1}, x, t ) < \epsilon.
		\end{align*}
		
		Continuing this process, at some $ k^{th} $ step, $ k \in \mathbb{N} $, there exist    $ n_k \geq k $ \& $ s_{n_k} \in F $ such that 	$ \mathfrak{P}_\theta ( s_{n_k}, x, t ) < \epsilon,  \quad \text{for  any}  ~ t > 0  $.
		
		From this we can conclude, 
		$$ \underset{n \rightarrow \infty}{\lim} \mathfrak{P}_\theta ( s_{n_k}, x, t ) = 0, \quad \text{for any} ~ t > 0. $$ 
		
		Thus any sequence in $F$ has a convergent subsequence in $F$. Hence $F$  is sequentially compact.  
		
	Conversely, suppose that $F$  is sequentially compact.  Let $ ( f, \alpha ) \in \mathfrak{F} \times [ 0, \infty ) $ be the pair with respect to which $ \mathfrak{P}_{\theta} $ satisfies $ ( \mathfrak{P}_{\theta} 2) $.
		
		First we claim that,
		\begin{equation} \label{12.4-4}
			\forall r > 0, ~\exists s_1, s_2, \cdots, s_n \in F ~ \& ~ t_1,t_2,\cdots, t_n \in ( 0, \infty ) ~ \text{such that} ~ F \subset \underset{\iota = 1, 2, \cdots, n }{ \cup} \mathscr{B} ( s_\iota, r, t_\iota ). 
		\end{equation}
		
		On the contradictory, suppose there exists $ r > 0 $ such that for any finite number of eements $ \{ s_1, s_2, \cdots, s_n \} \subset F $ and $ \{ t_1,t_2,\cdots, t_n  \} \subset ( 0, \infty )  $, we have 
		$$ F \nsubseteq  \underset{\iota = 1, 2, \cdots, n }{ \cup} \mathscr{B} ( s_\iota, r, t_\iota ).$$
		
		Let $ s_1 \in F $ be an arbitrary element. Then $ F \nsubseteq    \mathscr{B} ( s_1, r, t_1 ) $ for some $ t_1 > 0 $ i.e. there exists $ s_2 \in F $ such that $ \mathfrak{P}_\theta ( s_1, s_2, t_1 ) \geq r $.
		
		Again we have,   $ F \nsubseteq    \mathscr{B} ( s_1, r, t_1 ) \cup \mathscr{B} ( s_2, r, t_2 ) $ for some $ t_2 > 0 $ i.e. there exists $ s_3 \in F $ such that $ \mathfrak{P}_\theta ( s_i, s_2, t_i ) \geq r $ for $ i = 1, 2 $.
		
		Continuing this process, we can construct a sequence $ \{ s_n \} $ in $F$ such that $ \mathfrak{P}_\theta ( s_i, s_m, t_i ) \geq r $ for $ m\geq n $ and $ i = 1, 2, \cdots, n$.
		
		This shows that $ \{ s_n \} $  is not have a  Cauchy sequence and hence from Proposition \ref{conv-cauchy relation} does not exist  any convergent subsequence, which contradict the definition of sequentially compactness.

		Next suppose that $\{ \mathcal{O}_\iota\}_{\iota \in \Delta} $ be an arbitrary family of open subsets in $  \mathcal{X} $ such that 
		\begin{equation} \label{13.4-1}
			F \subset \underset{\iota \in \Delta }{ \cup} \mathcal{O}_\iota.
		\end{equation}
		
		We claim that,
		\begin{equation}\label{13.4-2}
			\exists r_0 > 0   :  \forall a\in F, ~ \forall t > 0, ~ \exists i \in \Delta ~\text{with} ~ \mathscr{B} ( a, r_0, t ) \subset \mathcal{O}_i. 
		\end{equation}
		
		On the contradictory, suppose that for any $ r> 0 $, there exist $ x_r \in F $ and $ t' > 0 $ such that
		\begin{equation}\label{13.4-3}
			\mathscr{B} ( s_r, r, t' ) \nsubseteq \mathcal{O}_i, \quad \forall i \in \Delta. 
		\end{equation}
		
		In particular, for all $ n \in \mathbb{N} $, there exists $ s_n \in F $ such that 
		\begin{equation}\label{13.4-4}
			\mathscr{B} ( s_n, \frac{1}{ n}, t' ) \nsubseteq \mathcal{O}_i, \quad \forall i \in \Delta. 
		\end{equation}
		
		Using sequentially compacteness, we can construct a subsequence $ \{ s_{n_k}\} $ of $ \{ s_n \} $ such that
		\begin{equation}  \label{13.4-5}
			\underset{k \rightarrow \infty}{\lim} \mathfrak{P}_\theta ( s_{n_k}, x, q ) = 0, \quad \text{for any} ~ q > 0
		\end{equation}
		
		where $ x $ is a certain element in $F$. Hence (\ref{13.4-1}) ensures the existence of some $ j \in \Delta $ such that $ x \in \mathcal{O}_j $.
		
		Since $ \mathcal{O}_j $ is an open set in $  \mathcal{X} $, so for $ t' > 0 $ there exist $ r_0 > 0 $   such that
		$$ \mathscr{B} ( x, r_0, t' ) \subset \mathcal{O}_j. $$
		
		In particular, for $ n_k \in \mathbb{N} $, take some $z \in \mathscr{B} ( s_{n_k}, \frac{1}{ n_k}, t' )$. So, $  \mathfrak{P}_\theta ( s_{n_k}, z, t' ) < \frac{1}{ n_k} $.
		
		Since, $ \mathfrak{P}_\theta ( s_{n_k}, z, u' ) $ is a continuous function ot $ u $, for all $ u > 0 $, so  there exists $ t'' \in ( 0, t' ) $ such that 
		$$   \mathfrak{P}_\theta ( s_{n_k}, z, t'' ) < \frac{1}{ n_k}. $$
		
		From (\ref{13.4-5}), in particular, we have 
		\begin{equation}  \label{13.4-6}
			\underset{k \rightarrow \infty}{\lim} \mathfrak{P}_\theta ( s_{n_k}, x, t' -t'' ) = 0.
		\end{equation}
		
		Therefore, using both ($\mathcal{B}3$) and (\ref{13.4-6}), we get
		$$ 	\underset{k \rightarrow \infty}{\lim} \theta \bigg( \mathfrak{P}_\theta ( s_{n_k}, x, t' -t'' ),  \mathfrak{P}_\theta ( s_{n_k}, z, t'' )  \bigg) < 	\underset{k \rightarrow \infty}{\lim} \theta \bigg( \mathfrak{P}_\theta ( x_{n_k}, x, t' -t'' ),  \frac{1}{ n_k}  \bigg) = 0.  $$
		
		Now, as $f$ satisfies $ (\mathscr{F}2)$, so the above inequality ensures the existence of a natural number $Q $ such that  
		$$  f \bigg( \theta \bigg( \mathfrak{P}_\theta ( s_{n_k}, x, t' -t'' ),  \mathfrak{P}_\theta ( s_{n_k}, z, t'' )  \bigg) \bigg) < f(r_0) - \alpha, \quad \forall n_k \geq Q.  $$
		
		Again, for the condition $ (\mathfrak{P}_\theta 2) $, $ \mathfrak{P}_\theta ( x, z, t' ) > 0 $ gives
		\begin{align*}
			& f (\mathfrak{P}_\theta ( x, z, t' )   ) \leq f \bigg( \theta \bigg( \mathfrak{P}_\theta ( s_{n_k}, x, t' -t'' ),  \mathfrak{P}_\theta ( s_{n_k}, z, t'' )  \bigg) \bigg) + \alpha \\
			or ~~ & f (\mathfrak{P}_\theta ( x, z, t' )   ) < f ( r_0 ), \quad \forall n_k \geq Q. 
		\end{align*} 
		
		This implies 
		\begin{align*}
			& \mathfrak{P}_\theta ( x, z, t' )   <   r_0, \quad (\text{using} ~   (\mathscr{F}1))  \\
			or ~~ & z \in \mathscr{B} ( x, r_0, t' ).
		\end{align*}
		
		Thus we have, $ \mathscr{B} ( s_{n_k}, \frac{1}{ n_k}, t' ) \subset  \mathscr{B} ( x, r_0, t' ) $, which implies $ \mathscr{B} ( s_{n_k}, \frac{1}{ n_k}, t' ) \subset \mathcal{O}_j $. This is a contradiction to (\ref{13.4-4}). Thus  (\ref{13.4-3}) does not hold.
		
		Again (\ref{12.4-4}) \& (\ref{13.4-2}) together implies that, for any  $ j = 1, 2, \cdots, n $, there exist $ i_j \in \Delta $ such that 
		$$ \mathscr{B} ( s_j, r_, t' ) \subset \mathcal{O}_{ i_j}. $$
		
		Hence, we can write that
		$$ F \subset  \underset{j = 1,  \cdots, n }{ \cup} \mathcal{O}_{ i_j}. $$
		
		Therefore $ F $ is compact.
	\end{proof}

  	\section{Suzuki's type fixed point theorem and its application}
  	
  	To facilitate our study, in the Section, we define a new version of Suzuki's contraction principle   in the setting of generalized $ \theta $-parametric metric space. We present two other theorems Suzuki-Banach  and Suzuki-Kannan type fixed point theorem  as particular case of the main  theorem which generalizes the concept of Banach and Kannan fixed point theorem.  Then apply it  in solving a fractional differential equations that arise in dynamic economic systems to analyze an economic growth model.

  	\subsection{Suzuki's type fixed point theorem }
  	
  	At this Section, we are ready to unveil the core theorems on  Suzuki's type contraction mapping.

  	\begin{dfn}
  		Let $ \mathfrak{T} $ be a self mapping on a generalized $ \theta $-parametric metric space $ ( \mathcal{X}, \mathfrak{P}_{\theta}, \theta ) $    and $ \psi : [ 0, 1) \to ( \frac{1}{2}, 1 ] $ be a mapping as defined in Theorem \ref{suzuki in m s}. $\mathfrak{T}$ is said to be Suzuki's type contraction mapping if for some  $ u \in [ 0, 1) $, $\mathfrak{T} $ satisfies the following condition:
  		\begin{equation}\label{ST1}
  			\psi ( u ) \mathfrak{P}_{\theta} ( x, \mathfrak{T} y, s ) \leq \mathfrak{P}_{\theta} ( x, y, s ) \implies \mathfrak{P}_{\theta} ( \mathfrak{T} x, \mathfrak{T} y, s ) \leq u M  (x, y, s )  
  		\end{equation}
  		for all  $ s> 0$ and for each pair $  ( x, y ) \in X \times X $ where 
  		$$ M  (x, y, t ) = \max \left\{ \mathfrak{P}_{\theta} ( x, y, t), \mathfrak{P}_{\theta} ( x, \mathfrak{T} x, t), \mathfrak{P}_{\theta} ( y, \mathfrak{T}y, t), \dfrac{( \mathfrak{P}_{\theta} ( x, \mathfrak{T}y, t) + \mathfrak{P}_{\theta} ( y, \mathfrak{T}x,  t) )}{2} \right\}. $$ 
  	\end{dfn}

  	We are ready to state a Suzuki's type fixed point theorem in the next theorem.

  	\begin{thm} 
  		\label{suzuki in new m s} \rm(Suzuki's type fixed point theorem \rm)
  		If $ \mathfrak{T} $ is a Suzuki's type contraction mapping on a complete generalized $ \theta $-parametric metric space $ ( \mathcal{X}, \mathfrak{P}_{\theta}, \theta ) $ and $ \mathfrak{P}_{\theta}  $ is  continuous in its first and second variable. 	Then $ \mathfrak{T} $ has an unique fixed point in $ \mathcal{X} $.
  	\end{thm}
  	\begin{proof} 
  		Since $ \psi ( u ) \leq 1 $ for all $ u \in [ 0, 1) $, so 
  		\begin{equation}\label{ST2}
  			\psi ( u ) \mathfrak{P}_{\theta} ( x, \mathfrak{T} x, s ) \leq \mathfrak{P}_{\theta} (x, \mathfrak{T} x, s )
  		\end{equation}
  		holds for any $ x \in X $ and $  s > 0 $. Using (\ref{ST1}), we have
  		\begin{align*}
  			\mathfrak{P}_{\theta} ( \mathfrak{T} x, \mathfrak{T}^2 x, t )  
  			\leq & ~ u M  (x, \mathfrak{T} x, s ) \\
  			\leq & ~ u \max \left\{ \mathfrak{P}_{\theta} ( x, \mathfrak{T}x, t), \mathfrak{P}_{\theta} ( x, \mathfrak{T} x, t), \mathfrak{P}_{\theta} ( \mathfrak{T} x, \mathfrak{T}^2 x, t), \dfrac{( \mathfrak{P}_{\theta} ( x, \mathfrak{T}^2 x, t) + \mathfrak{P}_{\theta} ( \mathfrak{T} x, \mathfrak{T} x,  t) )}{2} \right\} \\
  			= & ~ u \max \left\{ \mathfrak{P}_{\theta} ( x, \mathfrak{T} x, t),  \mathfrak{P}_{\theta} ( \mathfrak{T} x, \mathfrak{T}^2 x, t)  \right\} \\
  			= & ~ u   \mathfrak{P}_{\theta} ( x, \mathfrak{T} x, t)
  		\end{align*}
  		holds for all  $ t > 0 $.  Hence
  		\begin{equation}\label{ST3}
  			\mathfrak{P}_{\theta} ( \mathfrak{T} x, \mathfrak{T}^2 x, t ) \leq u   \mathfrak{P}_{\theta} ( x, \mathfrak{T} x, t) 
  		\end{equation}
  		holds for each  $~ x \in \mathcal{X} $ and    for all $  t > 0 $. Let $ w_0 \in \mathcal{X} $ be a fixed element. Let us consider an iterative sequence $ \{ w_\iota \} $ defined by 
  		$ w_{\iota +1} = \mathfrak{T}^n ( w_\iota), ~~ \iota = 0, 1, 2, \cdots. $ Using (\ref{ST3}), we have
  		$$ \mathfrak{P}_{\theta} ( w_{ \iota +1}, w_{\iota}, t   ) \leq u   \mathfrak{P}_{\theta} ( w_{\iota}, w_{\iota-1}, t   ) \leq \cdots \leq u^\iota   \mathfrak{P}_{\theta} ( w_{0}, \mathfrak{T} w_{0}, t   )  ~~~ \mbox{for all} ~ t > 0 ~ \&  ~ \iota \in \mathbb{N}. $$
  		This yields 
  		\begin{equation}\label{ST4}
  			\underset{\iota \to \infty }{lim} ~ \mathfrak{P}_{\theta} ( w_{\iota+1}, w_{\iota}, t   ) = 0 ~~~ \text{for all} ~ t > 0.	
  		\end{equation}
  		Suppose   the pair $ (f, \alpha) \in \mathfrak{F} \times [ 0, \infty  ) $ satisfy axiom $ (\mathfrak{P}_{\theta} 2) $.  Then we have
  		\begin{equation}\label{ST6}
  			f ( \mathfrak{P}_{\theta} ( w_{\iota+2}, w_{\iota}, t   ) ) \leq f \left(  \theta \left( \mathfrak{P}_{\theta} ( w_{\iota +2}, w_{\iota+1}, \frac{t}{2}     ), \mathfrak{P}_{\theta} ( w_{\iota+1}, w_{\iota}, \frac{t}{2} ) \right)  \right) + \alpha ~~~ \text{for all} ~ t > 0 ~ \& ~ \iota \in \mathbb{N}.
  		\end{equation}
  		Since $ \theta $ is continuous, using the condition ($ \mathcal{B} 1 $), and relation (\ref{ST4}), we obtain 
  		\begin{align*}
  			& \underset{\iota \to \infty }{lim} ~ 	 \theta \left( \mathfrak{P}_{\theta} ( w_{\iota+2}, w_{\iota+1}, \frac{t}{2}     ), \mathfrak{P}_{\theta} ( w_{\iota+1}, w_{\iota}, \frac{t}{2} ) \right) = \theta( 0, 0) = 0     \\
  			\implies & f \left(  \theta \left( \mathfrak{P}_{\theta} ( w_{\iota+2}, w_{\iota+1}, \frac{t}{2}     ), \mathfrak{P}_{\theta} ( w_{\iota+1}, w_{\iota}, \frac{t}{2} ) \right)  \right) + \alpha = - \infty   ~~~ (\text{using} ~~ (\mathfrak{F}2))  \\
  			\implies & \underset{\iota \to \infty }{lim} ~ f ( \mathfrak{P}_{\theta} ( w_{ \iota +2}, w_{\iota}, t   ) ) \leq - \infty ~~~ (\text{using} ~~ (\ref{ST6})) \\
  			\implies &  \underset{\iota \to \infty }{lim} ~ f ( \mathfrak{P}_{\theta} ( w_{ \iota +2}, w_{ \iota }, t   ) ) = - \infty \\
  			\implies & \underset{\iota \to \infty }{lim} ~ \mathfrak{P}_{\theta} ( w_{\iota+2}, w_{\iota}, t   ) = 0 ~~~ (\text{using} ~~ (\mathfrak{F}2)).  
  		\end{align*}
  		holds for all $ t > 0 $. 
  		Continuing in this process $p $ times,   we obtain 
  		\begin{equation}\label{ST5}
  			\underset{\iota \to \infty }{lim} ~ \mathfrak{P}_{\theta} ( w_{\iota+p}, w_{ \iota }, t   ) = 0 ~~~ \text{for all} ~ t > 0 ~ \& ~ p = 1, 2, \cdots.	
  		\end{equation}
  		This shows that $ \{ w_\iota\} $ is a Cauchy sequence in $X$ and hence it converges to some $z \in \mathcal{X} $. Again we have,
  		$$ \mathfrak{P}_{\theta} ( z, \mathfrak{T} z, t ) = \underset{\iota \to \infty }{lim} ~ \mathfrak{P}_{\theta} ( w_{\iota }, \mathfrak{T} w_{\iota}, t   ) = \underset{\iota \to \infty }{lim} ~ \mathfrak{P}_{\theta} ( w_{\iota }, w_{\iota+1}, t   ) = 0 ~~~ \text{for all} ~ t > 0. $$
  		This implies $ ~ \mathfrak{T} z = z ~$ i.e.  $ z $ is a fixed point for $ \mathfrak{T} $. Moreover, if $ \mathfrak{T} $ has another fixed point $ \eta \in \mathcal{X} $ then 
  		$$ \mathfrak{P}_{\theta} ( z,  \mathfrak{T}  z, t ) = 0 = \mathfrak{P}_{\theta} ( \eta, \mathfrak{T}  \eta, t ) ~~~ \text{holds for any} ~  t >   0. $$
  		Finally applying the realtion (\ref{ST1}), we  obtain
  		\begin{align*}
  			& 	\mathfrak{P}_{\theta} ( \mathfrak{T}  z, \mathfrak{T}  \eta, t ) \leq u \max \left\{ \mathfrak{P}_{\theta} ( z, \eta, t), \mathfrak{P}_{\theta} ( z, \mathfrak{T} z, t), \mathfrak{P}_{\theta} ( \eta, \mathfrak{T}  \eta, t), \dfrac{  \mathfrak{P}_{\theta} ( z, T \eta, t) + \mathfrak{P}_{\theta} ( \eta, \mathfrak{T}  z,  t)  }{2} \right\}  ~~~ t > 0 \\
  			\implies & \mathfrak{P}_{\theta} (   z,   \eta, t ) \leq u  \mathfrak{P}_{\theta} (   z,   \eta, t )  ~~~ t > 0 \\
  			\implies & \mathfrak{P}_{\theta} (   z,   \eta, t ) = 0 ~~~ t > 0 ~~~ (\text{since} ~ 0 \leq u < 1)\\
  			\implies & z = \alpha. 
  		\end{align*}
  		Hence the proof is complete.   
  	\end{proof}

  	The followings are Suzuki-Banach  and Suzuki-Kannan type fixed point results in generalized $ \theta $-paramatric metric space which can be derived as special cases of Theorem \ref{suzuki in new m s}.
  	
  	\begin{thm}\label{suzuki-banach in new m s} \rm(Suzuki-Banach type fixed point theorem\rm)
  		Let $ ( \mathcal{X}, \mathfrak{P}_{\theta}, \theta ) $ be a complete generalized $ \theta $-paramatric metric space  and $ \mathfrak{P}_{\theta}  $ be  continuous in its first and second variable.  
  		If there exists $ u \in [0, 1 ) $ with respect to which $ \mathfrak{T} $ satisfies   Suzuki-Banach type contraction condition 
  		\begin{equation}\label{SBT1}
  			\psi ( u ) \mathfrak{P}_{\theta} ( x, \mathfrak{T}  \eta, s ) \leq \mathfrak{P}_{\theta} ( x, \eta, s ) \implies \mathfrak{P}_{\theta} ( \mathfrak{T}  x, \mathfrak{T}  \eta, s ) \leq u \mathfrak{P}_{\theta}  (x, \eta, s )  
  		\end{equation}
  		for all  $ s> 0$ and for each pair $  ( x, \eta ) \in \mathcal{X} \times \mathcal{X} $ where $ \psi : [ 0, 1) \to ( \frac{1}{2}, 1 ] $ is the  function defined  in Theorem \ref{suzuki in m s}. Then $ \mathfrak{T} $ has an unique fixed point in $ \mathcal{X} $.	
  	\end{thm}

  	\begin{thm}\label{suzuki-kannan in new m s} \rm(Suzuki-Kannan type fixed point theorem\rm)
  		Let $ ( \mathcal{X}, \mathfrak{P}_{\theta}, \theta ) $ be a complete generalized $ \theta $-parametric metric space  and $ \mathfrak{P}_{\theta}  $ be  continuous in its first and second variable,  
  		If there exists $ u \in [0, 1 ) $ with respect to which $ \mathfrak{T} $ satisfies   Suzuki-Kannan type contraction condition 
  		\begin{equation}\label{SKT1}
  			\psi ( u ) \mathfrak{P}_{\theta} ( x, \mathfrak{T} \eta, s ) \leq \mathfrak{P}_{\theta} ( x, \eta, s ) \implies \mathfrak{P}_{\theta} ( \mathfrak{T} x, \mathfrak{T} \eta, s ) \leq \frac{u}{2} \left( \mathfrak{P}_{\theta}  (x, \mathfrak{T} \eta, s ) + \mathfrak{P}_{\theta} ( \eta, \mathfrak{T} x, s ) \right) 
  		\end{equation}
  		for all  $ s> 0$ and for each pair $  ( x, \eta ) \in \mathcal{X}^2   $ where $ \psi : [ 0, 1) \to ( \frac{1}{2}, 1 ] $ is the  function defined  in Theorem \ref{suzuki in m s}. Then $ \mathfrak{T} $ has an unique fixed point in $ \mathcal{X} $.	
  	\end{thm}

  	Now we provide an example to highlight our proven results.

  	\begin{eg}
  		Let $ \mathcal{X} = \{ ( 3, 3 ), ( 7, 3 ), ( 3, 7 ), ( 7, 9 ) \} .$ Consider the function $ B : \mathcal{X} \times \mathcal{X} \to [0, \infty ) $ by $ B ( \xi, y ) = ( \xi_1 - y_1 )^2 + ( \xi_2 - y_2 ) ^2 $ for all $ \xi = ( \xi_1, \xi_2 ), ~ y = ( y_1, y_2 )  \in \mathcal{X} .$ 
  		Then $ B $ is a $ b$-metric on $\mathcal{X} $ with constant coefficient $ K = 2 $ and hence $ \mathfrak{P}_{\theta} ( \xi, y, t ) = \dfrac{B ( \xi, y ) }{t} $ for all $ \xi, y \in \mathcal{X} $ and $ t > 0 $ is a generalized $ \theta $-parametric metric on $ \mathcal{X} $  (see Example \ref{ex via b-metric}). We define a mapping $ \mathfrak{T} : \mathcal{X} \to \mathcal{X} $ by $  \mathfrak{T} ( 3, 3 ) = ( 3, 3) = \mathfrak{T} ( 7, 3 ) = \mathfrak{T} ( 3, 7 ), ~ \mathfrak{T} ( 7, 9 ) = ( 7, 3). $
  		By routine calculation we can check that for $ r = \frac{7}{8} $ the relation  (\ref{ST1}) holds for all $ \xi, y \in \mathcal{X} $ and  $ t > 0 $.  Theorem \ref{suzuki in new m s} ensures the existence and uniqueness of a fixed point for $\mathfrak{T} $. This is clearly $ ( 3, 3 ) $.
  		In particular case, $ \mathfrak{P}_{\theta} $  satisfies all conditions of  Theorem \ref{suzuki-banach in new m s} and \ref{suzuki-kannan in new m s}. Further, if we consider $ Y = \mathcal{X} \cup \{ ( 9 , 7 ) \} $ and define a function on $Y$ by $ \mathfrak{S} ( y ) = \mathfrak{T} ( y) ~~ \text{if} ~ y \in \mathcal{X} ~~~ \text{and} ~~~ \mathfrak{S} ( 9 , 7 ) = ( 3, 7 ). $ It is not hard to see for $ x = ( 7, 9 ) $ and $ y = ( 9 , 7 ) $, $ ~   \mathfrak{P}_{\theta} ( \mathfrak{S} x, \mathfrak{S} y , t ) \leq r \mathfrak{P}_{\theta} ( x, y, t) $ holds for all $ \frac{7}{8} \leq r < 1 .$
  		One can observe that 
  		$ \psi ( r ) \mathfrak{P}_{\theta} ( x, \mathfrak{S} x, t ) \leq \mathfrak{P}_{\theta} ( x, y, t ) $ has no solution for $ r \in [ 0, 1 ) $. Thus $\mathfrak{S} $ in $Y$ satisfies non-Suzuki type condition. Hence the existence criterion of fixed point for $\mathfrak{S}$ can not be concluded from Theorem \ref{suzuki in new m s} even though $\mathfrak{S} $ has unique fixed point $ ( 3, 3 ) $ in $Y$. 
  	\end{eg}

  	\subsection{Application to  Economic Growth Modeling via Caputo Fractional Differential Approach }

  	Fractional differential equation  is applied to various fields of science and engineering. Specially the Caputo sense—offer a powerful tool for modeling economic growth by incorporating memory effects and dynamic flexibility.  Following   result  provides a rigorous mathematical foundation for analyzing such economic models. This application underscores the practical significance of our theoretical findings. By integrating generalized $\theta $-parametric metric spaces with fractional economic models, our results provide a powerful mathematical foundation for studying economic growth and also highlight the practical impact of fixed point theory and fractional differential equations in economics.

  	For, we consider a fractional differential equation
  	
  	\begin{equation}\label{app1}
  		D^\eta ( f ( \tau ) ) = g ( \tau, f ( \tau ) ), ~~ 0 < \tau < 1, ~~ 1 <  \eta \leq 2	
  	\end{equation}
  	
  	subject to  the  integral boundary conditions
  	\begin{equation}\label{app2}
  		f ( 0 ) = 0, ~ I f ( 1 ) = f' ( 0 ) 
  	\end{equation} 
  	
  	where $ D^\eta $ represents the Caputo fractional derivative of order $ \eta $ and defined by
  	\begin{equation}\label{app3}
  		D^\eta ( f ( t ) ) = \dfrac{1}{\Gamma ( i - \eta )} \int_{0}^{t} ( t- s)^{i-\eta - 1} f^i ( m ) dm 
  	\end{equation} 
  	
  	such that $ i -1 < \eta < i, ~ i = [\eta] + 1$ and $ f : [0, 1 ] \times \mathbb{R} \to [0, \infty ) $ is a continuous function and $ I^\eta f $ denotes the Reimann-Liouville fractional  integral of order $ \eta $  of a continuous function $ f : \mathbb{R}^+ \to \mathbb{R} $ given by
  	$ I^\eta f ( \tau ) = \dfrac{1}{\Gamma ( \eta )} \int_{0}^{ \tau } ( \tau - s)^{ \eta - 1} f  ( m ) dm. $ The variable $ f ( \tau ) $ could represent the economic indicator that characterizes economic health of a region. The non-linear function $ g ( \tau, f ( \tau ) ) $ contributes to various factor (viz. innovation, government spending, education system) of  an economic growth model. The non-local effects and memory in the economic system is represented by the fractional order $ \eta $. The integral boundary conditions $ f ( 0 ) = 0 $ and $ I f ( 1 ) = f' ( 0 ) $ represents a beginning of an economic activity or observation period and a connection between the accrued value   over the  period $ 0 $ to $1$ \& the rate of change of the economic variable  at the beginning of the observation period(please see \cite{15,16,17}). \\

  	The following theorem serves as a cornerstone,   providing a rigorous existence result for the considered fractional differential equation (\ref{app1}) within our newly introduced framework.

  	\begin{thm}
  		Let $ \mathcal{X} = C [ 0, 1 ] $ with sup metric $\ss _\infty (  \upphi, \phi  ) = {\underset{s\in [0,1]}{\sup} ~ | \upphi ( s ) - \phi ( s ) |} $ for all $ \upphi, \phi \in \mathcal{X}  $. Then $ ( \mathcal{X}, \mathfrak{P}_{\theta}, \theta ) $ is a complete generalized $ \theta $-parametric metric space   where
  		$ \mathfrak{P}_{\theta} ( \upphi, \phi, q ) = \dfrac{ \ss_\infty (  \upphi, \phi  )}{q} $ for all $ \upphi, \phi \in \mathcal{X}  ~ \& ~ q > 0  $
  		with respect to   $ f ( \xi ) = \ln \xi, ~ \forall \xi > 0 $, $ \alpha = 0$ and $ \theta ( \eth_1, \eth_2 ) = \eth_1 +  \eth_2 $ for all $ \eth_1, \eth_2 \in [ 0, \infty ) $.
  		
  		Consider the non-linear fractional differential equation (\ref{app3}) and suppose the following conditions are fulfilled:
  		\begin{enumerate}[(i)]
  			\item  $ f : [ 0, 1] \times \mathbb{R} \to \mathbb{R} $ is continuous;
  			\item if there exists $ r \in [ 0, 1 ) $ such that for every $ \Xi, \chi \in \mathcal{X}  $
  			\begin{equation}\label{app3*}
  				\psi ( r ) \mathfrak{P}_{\theta} ( \Xi, H \Xi, t ) \leq \mathfrak{P}_{\theta} ( \Xi, \chi, t ) 
  			\end{equation}
  			holds for  all $ t > 0 $ where $ H : \mathcal{X} \to \mathcal{X} $ is defined by 
  			\begin{equation}\label{app4}
  				H ( \Xi ( t ) ) = \dfrac{1}{\Gamma ( \eta )} \int_{0}^{t} ( t- s)^{ \eta - 1} f (s, \Xi (s)) ds +   \dfrac{2 t }{\Gamma ( \eta )} \int_{0}^{1}  \left( \int_{0}^{s} ( s - q )^{ \eta - 1} f  ( q, \Xi ( q ) ) dq    \right) ds.
  			\end{equation} 
  			The given function  $ f $ satisfy
  			\begin{equation}\label{app5}
  				|  f ( t, \chi_1 ( t ) ) - f ( t, \chi_2 ( t ) ) | \leq r \dfrac{\Gamma ( \eta + 1 ) }{4} | \chi_1 ( t ) - \chi_2 ( t ) | ~~ ~~ \text{for all } ~ t > 0. 
  			\end{equation} 
  		\end{enumerate}
  		Then the equation (\ref{app1}) has a unique solution in $\mathcal{X}$. 
  	\end{thm}
  	\begin{proof}
  		A solution $ h \in \mathcal{X} $ of (\ref{app1})   also satisfies  the integral equation (\ref{app4}). Suppose $\Xi, \chi \in \mathcal{X}  $ in such a way that the integral operator $H$ defined in (\ref{app4}) satisfies the condition (\ref{app3*}). For all $ t \in [ 0, 1 ] $
  		\begin{align*}
  			& | H ( \Xi ( t ) ) -  H ( \chi ( t ) ) | \\
  			= & ~  | \dfrac{1}{\Gamma ( \eta )}\int_{0}^{t} ( t- s)^{ \eta - 1} \left[  f (s, \Xi (s) ) -  f (s, \chi (s) )  \right] ds + \\
  			& \hskip 80 pt  \dfrac{2 t }{\Gamma ( \eta )} \int_{0}^{1} \left(  \int_{0}^{s} ( s - m)^{ \eta - 1} [ f  ( q, \Xi ( q ) ) -  f  ( q, \chi ( q ) ) ] dq \right) ds   |   \\
  			\leq & ~ \dfrac{1}{\Gamma ( \eta )} \int_{0}^{t} ( t- s)^{ \eta - 1}  |  f (s, \Xi (s) ) -  f (s, \chi (s) )  | ds + \\
  			&  \hskip 80 pt  \dfrac{2 t }{\Gamma ( \eta )} \int_{0}^{1} \left(  \int_{0}^{s} ( s - m)^{ \eta - 1} | f  ( q, \Xi ( q ) ) -  f  ( q, \chi ( q ) )  | dq   \right) ds \\
  			\leq & ~  \dfrac{1}{\Gamma ( \eta )} \int_{0}^{t} ( t- s)^{ \eta - 1}  r \dfrac{\Gamma ( \eta + 1 ) }{4} | \Xi ( s ) - \chi ( s ) | ds + \\
  			&  \hskip 80 pt  \dfrac{2 t }{\Gamma ( \eta )} \int_{0}^{1} \left(  \int_{0}^{s} ( s - q)^{ \eta - 1} r \dfrac{\Gamma ( \eta + 1 ) }{4}  |   \Xi ( q )   -    \chi ( q )    | dq   \right) ds \\
  			= & ~ \dfrac{r  \Gamma ( \eta + 1 ) }{ 4 \Gamma ( \eta )} \left[  \int_{0}^{t} ( t- s)^{ \eta - 1}  | \Xi ( s ) - \chi ( s ) | ds  + 2 t \int_{0}^{1} \left(  \int_{0}^{s} ( s - q )^{ \eta - 1}    |   \Xi ( q )   -    \chi ( q )    | dq   \right) ds \right]  \\
  			\leq & ~ \dfrac{r    \eta   }{ 4  }  \underset{u \in [ 0, 1]}{\sup}  |   \Xi ( u )   -    \chi ( u ) | \left( \int_{0}^{t} ( t- s)^{ \eta - 1} ds + 2t  \int_{0}^{1} \left(  \int_{0}^{s} ( s - q)^{ \eta - 1}  dq   \right) ds \right)  \\
  			\leq & ~ \dfrac{r   \eta  }{ 4  } ~ d_\infty ( \Xi, \chi ) ~ \left( \dfrac{1 + 2 t}{\eta}   \right)  
  			<   ~ r d_\infty ( \Xi, \chi ). 
  		\end{align*}
  		This implies \begin{align*}
  			& \underset{u \in [ 0, 1]}{\sup}  |   H( \Xi ( t ) )  -    H( \chi ( t )) | \leq  	r \ss_\infty ( \Xi, \chi ) \\&  \mbox{or} ~\dfrac{ \ss_\infty ( H (\Xi), H ( \chi) )}{l} \leq r ~ \dfrac{ \ss_\infty ( \Xi, \chi )}{l} ~~~~ \text{for all } ~ l > 0 \\& \mbox{ or}~  \mathfrak{P}_{\theta} (  H (\Xi), H ( \chi), l ) \leq r ~ \mathfrak{P}_{\theta} (  \Xi,   \chi, l ) ~~~~ \text{for all } ~ l > 0.
  		\end{align*}
  		
  		Thus $ H $ satisfies the Suzuki-Banach type contraction condition. By Theorem \ref{suzuki-banach in new m s}, we conclude the existence of an unique solution for (\ref{app1}). 
  	\end{proof}

  	\section*{Conclusion:}
  	In this article, we introduce the concept of a generalized $ \theta$-parametric metric space, which extends and unifies the notions of both $ \theta$-metric spaces and parametric metric spaces. We undertake a detailed study of its fundamental properties, including convergence and Cauchy sequences, and provide non-trivial examples that illustrate and support our findings. 
  	
  	The establishment of a Suzuki-type fixed point theorem within this new framework is a  major highlight of our work. As a compelling application, we employ this result to analyze an economic growth model by investigating the existence of solutions to Caputo fractional differential equations which   demonstrates the practical significance of our theoretical contributions.
  	 This work not only broadens the scope of parametric metric spaces but also bridges abstract mathematical theory with real-world applications.

   	The general definition of a generalized $\theta$-parametric metric space,  while broad and unifying, does not in itself guarantee standard topological properties such as openness of open balls \& closedness  of closed balls or continuity of the distance function. These phenomenons have been demonstrated through explicit counterexamples. This suggests that additional assumptions (e.g., continuity in the parameter, or properties of the $\mathcal{B}$-action) may be necessary to recover such features. Identifying minimal sufficient conditions remains an interesting and relevant direction for future exploration in this direction. 
  	~\\

  	\noindent
  	\textbf{Data Availability Statement:} 
  	Data sharing not applicable to this article as no datasets were generated or analysed during the current study.\\
  	
  	\noindent
  	\textbf{Author Contributions:}
  	Conceptualization, Investigation, Methodology, Writing-original draft Abhishikta Das; 
  	Conceptualization,  Methodology,  Writing–review \& editing  Hemanta Kalita;
  	Conceptualization,  Funding, Validation, Writing–review \& editing  Mohammad Sajid; 
  	Conceptualization,  Methodology,   Writing–review \& editing   T. Bag.  \\
  	
  	\noindent
  	\textbf{Acknowledgement:}
  	The Researchers would like to thank the Deanship of Graduate Studies and Scientific Research at Qassim University for financial support (QU-APC-2025).  \\
  	
  	\noindent 
  	\textbf{Conflict of interest:}  The authors declare there is no conflicts of interest.

  \end{document}